\documentclass[10pt,oneside]{amsart}
\usepackage{tikz}
\usetikzlibrary{matrix,arrows,decorations.pathmorphing}
\usepackage{amssymb,mathtools,booktabs,enumitem}
\usepackage[letterpaper,textwidth=360pt,textheight=584pt,top=116pt,headheight=12pt,headsep=12pt]{geometry}
\usepackage[unicode,hidelinks]{hyperref}
\allowdisplaybreaks[2]
\AtBeginDocument{%
 \setlength{\oddsidemargin}{\dimexpr(\paperwidth-\textwidth)/2-1in\relax}%
 \setlength{\evensidemargin}{\oddsidemargin}%
}
\newcommand{\uq}{\mathbf{\overline{u}}_q(\mathfrak{sl_2})}
\newcommand{\kk}{\Bbbk}
\newcommand{\Z}{\mathbb Z}
\newcommand{\res}[1]{\langle #1\rangle_p}
\newcommand{\cop}{\mathrm{cop}}
\newcommand{\Span}{\operatorname{span}}

\newtheorem{proposition}{Proposition}[section]
\newtheorem{lemma}[proposition]{Lemma}
\newtheorem{corollary}[proposition]{Corollary}
\newtheorem{theorem}[proposition]{Theorem}

\theoremstyle{definition}

\newcommand{\selabel}[1]{\label{se:#1}}

\def\<{\leqslant}
\def\>{\geqslant}
\def\a{\alpha}
\def\b{\beta}

\def\g{\gamma}

\def\e{\varepsilon}

\def\ot{\otimes}

\date{}

\begin{document}
\title[Representations of the quantum double]{Representations of the quantum double of the restricted quantum group $\uq$}
\author{Hua Sun$^1$}
\address{College of Mathematical Science, Yangzhou University,
	Yangzhou 225002, China}
\email{huasun@yzu.edu.cn}
\author{Huixiang Chen$^1$}
\address{College of Mathematical Science, Yangzhou University,
	Yangzhou 225002, China}
\email{hxchen@yzu.edu.cn}
\author{Xiaoming Li$^2$}
\address{School of Statistics and Data Science, Taizhou University,
	Taizhou 225300, China}
\email{xiaomingmath@126.com}
\author{Huan Jia$^3$}
\address{Department of Mathematics, Suqian University, Suqian 223800, China}
\email{huan.jia@squ.edu.cn}

\thanks{Corresponding author: Hua Sun; huasun@yzu.edu.cn}

\thanks{1. College of Mathematical Science, Yangzhou University,
	Yangzhou 225002, China}

\subjclass[2010]{16T99, 16E99}
%\thanks{2020 {\it Mathematics Subject Classification}. 16T05, 16G30 }
\keywords{Quantum double, restricted quantum group, irreducible representation, highest weight module}
\begin{abstract}
Let $q$ be a primitive $2p$-th root of unity with $p\ge2$, and let $H=D(\uq)$ be the Drinfeld double of the restricted quantum group $\uq$. We describe the dual Hopf algebra and the cross relations of $H$. A change of generator gives a classification of the simple $H$-modules.  Their isomorphism classes are indexed by $\mathbb Z_{2p}^{\,2}$. We construct their projective covers, compute their dimensions, and determine which simple modules are projective. Finally, we prove that $H$ is strictly wild for every $p\ge2$ by constructing a fully faithful exact functor from finite-dimensional $\kk\langle x,y\rangle$-modules to finite-dimensional $H$-modules.
\end{abstract}
\maketitle

\section{Introduction and preliminaries}\selabel{1}
Finite-dimensional quantum groups at roots of unity give examples of nonsemisimple tensor categories. Their representation theory involves both the structure of individual modules and the behavior of modules under tensor products. The restricted quantum group $\uq$ is one of the basic examples in this setting. Its modules were studied by Suter and Xiao \cite{Suter,Xiao}, while tensor product decompositions and representation rings were investigated in \cite{KonSai,Su Yang}. These results make $\uq$ a useful starting point for studying how constructions on Hopf algebras affect simple modules, projective covers, and representation type.

The Drinfeld double is particularly relevant to this question. For a finite-dimensional Hopf algebra $A$, the double $D(A)$ is quasitriangular, and its representations describe compatible $A$-module and $A$-comodule structures through Yetter--Drinfeld modules \cite{Ka,Mo}. Categorically, these representations realize the Drinfeld center of the category of finite-dimensional $A$-modules. Thus the double relates explicit Hopf algebra calculations to the study of braided finite tensor categories and their centers. Determining its simple and projective modules gives concrete information about the objects and extensions in such a center, while its representation type measures the complexity of classifying those objects. The double of $\uq$ offers a setting in which these questions can be examined through generators, relations, and explicit module constructions.

There is also a connection with the algebraic study of logarithmic conformal field theory. Restricted quantum groups and their quasi-Hopf modifications occur in the representation theory of the triplet vertex operator algebra $W(p)$; see \cite{B.A,K.N.A,CGR}. In this context, nonsemisimple module categories, projective objects, and braided tensor structures are central issues. Studying $D(\uq)$ contributes explicit examples to the related theory of finite-dimensional quasitriangular Hopf algebras. It also allows one to compare two aspects of representation theory that need not have the same complexity: the classification of simple and projective modules, and the classification of arbitrary indecomposable modules.

The representation theory of doubles has been studied for Taft algebras \cite{H.C.X99,Ch2,Ch3,Ch4}, generalized Taft algebras \cite{EGST2006,EGST2019}, and pointed rank one Hopf algebras \cite{S.C.Z,S.C.z}. In this paper we consider $H=D(\uq)$ over an algebraically closed field $\kk$ of characteristic zero, with $q$ of order $2p$ and $p\ge2$. We classify the simple modules, construct their projective covers and compute their dimensions, and prove that $H$ is strictly wild for every $p\ge2$.

We first describe the dual Hopf algebra and derive the cross relations of $H$. For the representation theory, the useful generators are
\[
 X=E+\frac{\beta\gamma^2K}{q-q^{-1}},\qquad
 Y=F-\frac{q^2\alpha K^{-1}}{q-q^{-1}}.
\]
Together with $K,\alpha,\beta,\gamma$, these generators give an ordered basis of $H$. We use this basis to construct $p^2$-dimensional induced modules $Z(i,j)$ and identify their unique maximal submodules. The corresponding quotients $L(i,j)$ have rectangular bases. Our first main result is the following classification, proved in Theorem~\ref{thm:classification}. Write $\res t$ for the least nonnegative residue of $t$ modulo $p$.

\newtheorem*{introresultA}{Theorem A}
\newtheorem*{introresultB}{Theorem B}
\newtheorem*{introresultC}{Theorem C}
\begin{introresultA}
The simple $H$-modules are the modules $L(i,j)$, indexed by $(i,j)\in\mathbb Z_{2p}^{\,2}$. Their dimensions are
\[
 \dim L(i,j)=
 \begin{cases}
 p^2,&i\text{ odd},\\
 \bigl(1+\res{j+i/2}\bigr)\bigl(1+\res{i/2-j}\bigr),&i\text{ even}.
 \end{cases}
\]
In particular, there are exactly $4p^2$ simple-module isomorphism classes.
\end{introresultA}

For $\lambda=(i,j)\in\Lambda=\mathbb Z_{2p}^{\,2}$, let $L_\lambda=L(i,j)$ and $d_\lambda=\dim L_\lambda$. The construction in Section~\ref{sec:projectives} lifts a rank-one idempotent in the $\lambda$-factor of $H/J(H)$ to the idempotent $e_\lambda$ in \eqref{eq:projective-definition}.

\begin{introresultB}
The modules $P_\lambda=He_\lambda$ are the projective covers of the $L_\lambda$ and form a complete set of pairwise nonisomorphic indecomposable projective $H$-modules. Moreover,
\[
 {}_HH\cong\bigoplus_{\lambda\in\Lambda}P_\lambda^{\,d_\lambda},
 \qquad
 \dim P_\lambda=p^2\varepsilon(r(i,j))\varepsilon(s(i,j)),
\]
where $r(i,j)$ and $s(i,j)$ are the side lengths of the rectangular basis of $L(i,j)$, and $\varepsilon(t)=1$ if $t=p$ and $\varepsilon(t)=2$ if $t<p$. In particular, $L(i,j)$ is projective precisely when $r(i,j)=s(i,j)=p$.
\end{introresultB}

We encode a pair of linear operators in the generator actions of an $H$-modules. Theorem~\ref{thm:strict-wildness} gives the resulting embedding.

\begin{introresultC}
Let $R=\kk\langle x,y\rangle$. There is an $H$--$R$-bimodule $\mathcal M$, free of finite rank as a right $R$-module, such that
\[
 \mathcal M\otimes_R-:
 R\text{-}\mathrm{mod}\longrightarrow H\text{-}\mathrm{mod}
\]
is fully faithful and exact on finite-dimensional modules. The right rank is $3p-1$ for $p\ge3$ and $3$ for $p=2$. Consequently, $H$ is strictly wild for every $p\ge2$.
\end{introresultC}

The paper is organized as follows. Section~\ref{se:1} fixes notation and recalls the restricted quantum group. Section~\ref{se:3} describes the dual structure and the cross relations. Section~\ref{se:4} classifies simple modules, Section~\ref{sec:projectives} constructs their projective covers and computes their dimensions, and Section~\ref{sec:wildness} proves strict wildness.

Throughout, we work over an algebraically closed field $\Bbbk$ with char$(\Bbbk)$=0. Unless
otherwise stated, all algebras and Hopf algebras are
defined over $\Bbbk$, our references for basic concepts and notations about  Hopf algebras
are \cite{Ka, Mo}.
In particular, for a Hopf algebra, we will use $\e$, $\Delta$ and $S$ to denote the counit,
comultiplication and antipode, respectively.

Let $0\not=q\in \Bbbk$. For any nonnegative integer $n$, define $(n)_q$ by $(0)_q=0$ and $(n)_q=1+q+\cdots +q^{n-1}$ for $n>0$.
Observe that $(n)_q=n$ when $q=1$, and
$$
(n)_q=\frac{q^n-1}{q-1}
$$
when $q\not= 1$.
Define the $q$-factorial of $n$ by
$(0)!_q=1$ and
$(n)!_q=(n)_q(n-1)_q\cdots (1)_q$ for $n>0$.

The  $q$-binomial coefficients
$\binom{n}{i}_q$ for $0\<i\<n$
are defined by
$$\binom{n}{0}_q
=\binom{n}{n}_q
=1\ {\rm and} \
\binom{n}{i}_q
=\frac{(n)!_q}{(i)!_q(n-i)!_q}
$$
when
$(n-1)!_q\not = 0$ and $0<i<n$
(see \cite[Page 74]{Ka}).

For any integer $n$ and $q\in \Bbbk$ with $q\neq \pm 1$, set
$$
[n]=\frac{q^n-q^{-n}}{q-q^{-1}}.
$$
For a set $X$, $\sharp X$ denotes the number of elements in $X$.

Let $q\in \Bbbk$ with $q\neq\pm 1$.
The quantum enveloping algebra  $U_q=U_q(\mathfrak{sl}_2)$ of Lie algebra $\mathfrak{sl}_2$
is generated, as an algebra, by $E$, $F$, $K$ and $K^{-1}$ subject to the following relations:
$$\begin{array}{c}
KEK^{-1}=q^2E, \ KFK^{-1}=q^{-2}F, \ [E,F]=\frac{K-K^{-1}}{q-q^{-1}}, \ KK^{-1}=K^{-1}K=1.
\end{array}$$
$U_q$ is a Hopf algebra with the coalgebra structure and the antipode given by
$$\begin{array}{c}
\bigtriangleup(K)=K\otimes K, \ \bigtriangleup(E)=E\otimes K+1\otimes E, \ \bigtriangleup(F)=F\otimes 1+K^{-1}\otimes F,\\
\varepsilon(K)=1, \ \varepsilon(E)=\varepsilon(F)=0, \ S(K)=K^{-1}, \ S(E)=-EK^{-1}, \ S(F)=-KF.\\
\end{array}$$
Note that $K$ and $K^{-1}$ are group-like elements and $E,\ F$ are skew-primitive elements.
Now let $p\>2$ be an integer and let $q\in \Bbbk$ be a root of unity with the order $|q|=2p$. The ideal $I$ of $U_q$ generated by $E^p$, $F^p$ and $K^{2p}-1$
is a Hopf ideal of $U_q$, and hence one can form a quotient Hopf algebra
$\mathbf{\overline{u}}_q(\mathfrak{sl_2})=U_q/I$, see \cite{Suter,Xiao}. We call $\mathbf{\overline{u}}_q(\mathfrak{sl_2})$ the restricted quantum group. It is easy to see  dim$\mathbf{\overline{u}}_q(\mathfrak{sl_2})=2p^3$, and $\{E^iF^jK^l|0\<i,j\<p-1,l\in \mathbb{Z}_{2p}\}$ is a $\Bbbk$-basis of $\mathbf{\overline{u}}_q(\mathfrak{sl_2})$.

\section{Quantum double $D(\mathbf{\overline{u}}_q(\mathfrak{sl_2}))$}\selabel{3}

In this section, we describe the structure of the quantum double $D(\mathbf{\overline{u}}_q(\mathfrak{sl_2}))$.

Let $H$ be a finite dimensional Hopf algebra. Then its dual Hopf algebra $H^*$ is an $H$-bimodule with the left action $\rightharpoonup$ and right action $\leftharpoonup$ given by
$$(h\rightharpoonup f)(x)=f(xh) \  {\rm and} \  (f\leftharpoonup h)(x)=f(hx)$$ for all $h,x\in H$ and $f\in H^*$. Moreover, $H^*$ is a left(right) $H$-module algebra under the action $\rightharpoonup$ $(\leftharpoonup)$.

The Drinfeld double $D(H)$ of a finite dimensional Hopf algebra $H$ is a bicrossed product $H^{*cop}\bowtie H$. As coalgebra, $H^{*cop}\bowtie H=H^{*cop}\ot H$, the tensor product of coalgebra of $H^{*cop}$ and $H$. Denote by $f\bowtie h$ the element of $f\ot h$, $f\in H^{* cop}$, $h\in H$. Then the multiplication of $D(H)$ is given by $$(f\bowtie h)(g\bowtie x)=\sum f(h_1\rightharpoonup g\leftharpoonup S^{-1}(h_3))\bowtie h_2 x,$$ where $f,g\in H^{*cop}$ and $h,x\in H$, $S$ is the antipotde of $H$, see \cite{Ka,Mo}. $H^{*cop}$ and $H$ are Hopf subalgebra of $D(H)$ via the embeddings $H^{* cop}\hookrightarrow D(H)$, $f\mapsto f\bowtie 1$ and $H\hookrightarrow D(H)$, $h\mapsto \varepsilon\bowtie h$, respectively. Hence $D(H)$ is generated, as an algebra, by its two Hopf algebras $H^{*cop}$ and $H$ with the relations:
$$hf=\sum(h_1\rightharpoonup f\leftharpoonup S^{-1}(h_3))h_2, h\in H, f\in H^{*cop},$$ where $S$ is the antipode of $H$.

Now let $H=\mathbf{\overline{u}}_q(\mathfrak{sl_2})$. Then
the quantum double $D(\mathbf{\overline{u}}_q(\mathfrak{sl_2}))$  is the bicrossed product of $\mathbf{\overline{u}}_q(\mathfrak{sl_2})$ and of $\mathbf{\overline{u}}_q(\mathfrak{sl_2})^ {*cop}$.
$$D(\mathbf{\overline{u}}_q(\mathfrak{sl_2}))= \mathbf{\overline{u}}_q(\mathfrak{sl_2})^ {*cop}\bowtie \mathbf{\overline{u}}_q(\mathfrak{sl_2}).$$
Hence $\mathbf{\overline{u}}_q(\mathfrak{sl_2})$ and $\mathbf{\overline{u}}_q(\mathfrak{sl_2})^ {*cop}$
are Hopf subalgebras of $D(\mathbf{\overline{u}}_q(\mathfrak{sl_2}))$  via the identifications $y=\varepsilon\bowtie y$, $y\in \mathbf{\overline{u}}_q(\mathfrak{sl_2})$ and $f=f\bowtie 1$, $f\in \mathbf{\overline{u}}_q(\mathfrak{sl_2})^*$, respectively, where $\varepsilon$ is the unit element of $\mathbf{\overline{u}}_q(\mathfrak{sl_2})^*$. Let
$\{\overline{E^iF^jK^l}|0\<i,j\<p-1,l\in \mathbb{Z}_{2p}\}$ be the basis of $\mathbf{\overline{u}}_q(\mathfrak{sl_2})^*$ dual to the basis $\{E^iF^jK^l|0\<i,j\<p-1,l\in \mathbb{Z}_{2p}\}$.

\begin{lemma}\label{2.1} Let $0\<i,i'j,j'\<p-1$ and $l,l'\in \mathbb{Z}_{2p}$.
Then
$$\overline{E^iF^jK^l}\ \overline{ E^{i'}F^{j'}K^{l'}}=
q^{-2ij'}\binom{i+i'}{i'}_{q^2}\binom{j+j'}{j}_{q^2}\overline{E^{i+i'}F^{j+j'}K^{l+j'}}$$

if $i+i'<p$, $j+j'<p$ and $l'=l+i+j'$ in $\mathbb{Z}_{2p}$, and
$$\overline{E^iF^jK^l}\ \overline{ E^{i'}F^{j'}K^{l'}}=0$$ for other cases.
\end{lemma}
\begin{proof} Assume that $\overline{E^iF^jK^l}\ \overline{ E^{i'}F^{j'}K^{l'}}=\sum_{r,s=0}^{p-1}\sum_{t\in\mathbb{Z}_{2p}} \theta_{r,s,t}\overline{E^rF^sK^t}.$
For any $0\<r,s<p$ and $t\in \mathbb{Z}_{2p}$, a straightforward verification shows that
$$\begin{array}{rl}
&\Delta(E^rF^sK^t)\\
=&\sum_{x=0}^r\sum_{y=0}^sq^{2(x-r)(s-y)}\binom{r}{x}_{q^2}
\binom{s}{y}_{q^2}E^{r-x}F^{y}K^{y-s+t}\ot E^{x}F^{s-y}K^{r-x+t}.
\end{array}$$
Hence we have $$\begin{array}{rl}
&\theta_{r,s,t}\\
=&(\overline{E^iF^jK^l}\ \overline{ E^{i'}F^{j'}K^{l'}})(E^rF^sK^t)\\
=&\sum_{x=0}^r\sum_{y=0}^sq^{2(x-r)(s-y)}\binom{r}{x}_{q^2}\binom{s}{y}_{q^2}\delta_{i,r-x}\delta_{j,y}
\delta_{l,y-s+t}\delta_{i',x}\delta_{j',s-y}\delta_{l',r-x+t}.\\
=&\left\{
\begin{array}{ll}
\vspace{0.2cm}
q^{-2ij'}\binom{i+i'}{i'}_{q^2}\binom{j+j'}{j}_{q^2}& \begin{array}{l}\text{if }r=i+i',\ s=j+j',\\ t=l'-i=l+j'\text{ in }\mathbb Z_{2p},\end{array}\\
0, & {\rm otherwise}.\\
\end{array}\right.
\end{array}$$
Then, the lemma follows.
\end{proof}

The identity of the dual algebra $\uq^*$ is
 $\varepsilon=\sum_{i\in\mathbb Z_{2p}}\overline{K^i}.$
Put
\[
 \alpha=\sum_{i\in\mathbb Z_{2p}}\overline{EK^i},\qquad
 \beta=\sum_{i\in\mathbb Z_{2p}}\overline{FK^i},\qquad
 \gamma=\sum_{i\in\mathbb Z_{2p}}q^i\overline{K^i}.
\]

\begin{lemma}\label{2.2}
In $\mathbf{\overline{u}}_q(\mathfrak{sl_2})^*$, we have
$$\a^{p}=\b^p=0, \g^{2p}=\varepsilon,\ \a\b=q^{-2}\b\a,\ \a\g=q\g\a,\ \b\g=q\g\b.$$
\end{lemma}
\begin{proof}
We claim that $\a^i=(i)!_{q^2}\sum_{l\in \mathbb{Z}_{2p}}\overline{E^iK^l}$ for $0\<i<p$. It is trivial for $i=0$ and $i=1$. Now let $1\<i<p-1$. Then by the induction hypothesis and Lemma \ref{2.1}, we have
$$\begin{array}{rl}
\a^{i+1}&=(i)!_{q^2}\sum_{l,l'\in \mathbb{Z}_{2p}}\overline{E^iK^l} \ \overline{EK^{l'}}\\
&=(i)!_{q^2}\sum_{l\in \mathbb{Z}_{2p}}\binom{i+1}{1}_{q^2}\overline{E^{i+1}K^l}\\
&=(i+1)!_{q^2}\sum_{l\in \mathbb{Z}_{2p}}\overline{E^{i+1}K^l}.\\
\end{array}$$
Thus, we have proven the claim. Hence again by Lemma \ref{2.1},
$$\a^p=(p-1)!_{q^2}\sum_{l,l'\in \mathbb{Z}_{2p}}\overline{E^{p-1}K^l} \ \overline{EK^{l'}}=0.$$
Similarly, one can show that $\b^p=0$ and $\g^{2p}=\varepsilon.$ By Lemma \ref{2.1}, it is easy to check that $\a\b=q^{-2}\b\a$, $\a\g=q\g\a$ and $\b\g=q\g\b$.
\end{proof}

\begin{lemma}\label{2.3}
$\mathbf{\overline{u}}_q(\mathfrak{sl_2})^*$ is generated, as an algebra, by $\a$, $\b$ and $\g$.
\end{lemma}
\begin{proof} Let $0\<i,j<p$ and $k\in \mathbb{Z}_{2p}.$
By Lemma \ref{2.1} an argument similar to the proof of Lemma \ref{2.2} shows that
$$\a^i\b^j\g^k=(i)!_{q^2}(j)!_{q^2}q^{ki-2ij}\sum_{l\in\mathbb{Z}_{2p}}q^{kl}\overline{E^iF^jK^l}.$$
For any fixed $i$ and $j$, let $$A=\left(
       \begin{array}{cccc}
         1 & 1 & \cdots & 1 \\
         1 & q & \cdots & q^{2p-1} \\
         \vdots & \vdots & \ddots & \vdots \\
         1 & q^{2p-1} & \cdots & q^{(2p-1)(2p-1)} \\
       \end{array}
     \right)
.$$ Since ${\rm det}(A)\neq 0$, $A$ is an invertible matrix. This implies that  ${\rm span}\{\a^i\b^j\g^l|l\in \mathbb{Z}_{2p}\}={\rm span}\{\overline{E^iF^jK^l}|l\in \mathbb{Z}_{2p}\}$. Thus, the lemma follows.
\end{proof}

\begin{corollary}\label{2.4}
$\mathbf{\overline{u}}_q(\mathfrak{sl_2})^*$ has a $\Bbbk$-basis $\{\a^i\b^j\g^k|0\<i,j\<p-1, k\in \mathbb{Z}_{2p}\}$.
\end{corollary}

\begin{proof}
It follows from Lemmas \ref{2.2} and \ref{2.3}.
\end{proof}

For $0\le i,j\le p-1$, the normal-ordering identities in \cite[Lemma VI.1.3]{Ka} give
\[
 FE^i=E^iF-[i]E^{i-1}T_i,\qquad
 F^jE=EF^j-[j]F^{j-1}T'_j,
\]
where a correction term with zero index is omitted, and
\[
 T_i=\frac{q^{i-1}K-q^{1-i}K^{-1}}{q-q^{-1}},\qquad
 T'_j=\frac{q^{1-j}K-q^{j-1}K^{-1}}{q-q^{-1}}.
\]

\begin{lemma}\label{gs1}
If $1\<j\<i\<p-1$, then
$$F^jE^i=\sum_{t=0}^jE^{i-t}F^{j-t}f^{i,j}_t(K),$$ where $f^{i,j}_t(K)\in \Bbbk[ K,K^{-1}]$ for $0\<t\<j$ with $f^{i,j}_0(K)=1$ and $f^{i,j}_j(K)=(-1)^j[i][i-1]\cdots[i-j+1]T_{i-j+1}T_{i-j+2}\cdots T_i$.
\end{lemma}
\begin{proof}
 For $j=1$, it follows from the  above discussion. Now let $1\<j<i$ and assume that $F^jE^i=\sum_{t=0}^jE^{i-t}F^{j-t}f^{i,j}_t(K),$ where $f^{i,j}_t(K)\in \Bbbk[K,K^{-1}]$, $0\<t\<j$, with $f^{i,j}_0(K)=1$ and  $f^{i,j}_j(K)=(-1)^j[i][i-1]\cdots[i-j+1]T_{i-j+1}T_{i-j+2}\cdots T_i$. Then
$$\begin{array}{rl}
&F^{j+1}E^i\\
=&\sum_{t=0}^jFE^{i-t}F^{j-t}f^{i,j}_t(K)\\
=&\sum_{t=0}^j(E^{i-t}F-[i-t]E^{i-t-1}T_{i-t})F^{j-t}f^{i,j}_t(K)\\
=&\sum_{t=0}^jE^{i-t}F^{j-t+1}f^{i,j}_t(K)-\sum_{t=0}^j[i-t]E^{i-t-1}T_{i-t}F^{j-t}f^{i,j}_t(K)\\
=&\sum_{t=0}^jE^{i-t}F^{j-t+1}f^{i,j}_t(K)-\sum_{t=0}^j[i-t]E^{i-t-1}F^{j-t}T_{i+t-2j}f^{i,j}_t(K)\\
=&E^iF^{j+1}+\sum_{t=1}^jE^{i-t}F^{j-t+1}(f^{i,j}_{t}(K)-[i-t+1]T_{i+t-1-2j}f^{i,j}_{t-1}(K))\\
&-[i-j]E^{i-j-1}T_{i-j}f^{i,j}_{j}(K).\\
\end{array}$$
Hence $F^{j+1}E^i=\sum_{t=0}^{j+1}E^{i-t}F^{j+1-t}f^{i,j+1}_t(K)$, where $f^{i,j+1}_0(K)=1$,
$$f^{i,j+1}_t(K)=f^{i,j}_t(K)-[i-t+1]T_{i+t-1-2j}f^{i,j}_{t-1}(K)\in \Bbbk[K,K^{-1}]$$ for $1\<t\<j$, and $f^{i,j+1}_{j+1}(K)=-[i-j]T_{i-j}f^{i,j}_j(K)=(-1)^{j+1}[i][i-1]\cdots[i-j]T_{i-j}T_{i-j+1}\cdots T_i$.

 This completes the proof.
\end{proof}

\begin{lemma}\label{gs2}
If $0\<i<j\<p-1$, then
$$F^jE^i=\sum_{t=0}^iE^{i-t}F^{j-t}g^{i,j}_t(K),$$ where $g^{i,j}_t(K)\in \Bbbk[K,K^{-1}]$ for $0\<t\<i$ with $g^{i,j}_0(K)=1$ and $g^{i,j}_i(K)=(-1)^i[j][j-1]\cdots[j-i+1]T'_{j-i+1}T'_{j-i+2}\cdots T'_{j-2i+2}$ if $i>0$.
\end{lemma}
\begin{proof}
It is similar to Lemma \ref{gs1}
\end{proof}

\begin{proposition}\label{2.7}\label{2.8}\label{2.9}
In $\uq^*$, the following identities hold:
\begin{enumerate}[label=\textup{(\alph*)}]
\item[(a)]
$\Delta(\a)=\sum_{j=0}^{p-2}(-1)^jq^{-j^2}\b^j\g^{2(j+1)}\ot \a^{j+1}+\a\ot \varepsilon.$
\item[(b)]
$\Delta(\b)=\g^{-2}\ot \b+\b\ot \epsilon+(q^3+q)\b\ot \a\b+\b^2\g^2\ot (q^2\a^2\b+q^{-1}\a).$
\item[(c)]
$\Delta(\g)=\sum_{j=0}^{p-1}(-1)^jq^{-j^2}\b^j\g^{2j+1}\ot \a^{j}\g.$
\end{enumerate}
\end{proposition}
\begin{proof}
\textup{(a)}
By Lemma \ref{gs1}, $f^{j,j}_j(K)=\sum_{s=0}^jp^j_s K^{j-2s}$ for some $p^j_0,p^j_0,\cdots,p^j_j\in \Bbbk$, and $\sum_{s=0}^jp^j_{s}=0$, where $1\<j\<p-1.$ Similarly, $f^{j+1,j}_j(K)=\sum_{s=0}^j\varphi^j_sK^{j-2s}$ for some $\varphi^j_s\in \Bbbk$, where $1\<j\<p-2$. It is easy to check that $\sum_{s=0}^j\varphi^j_s=(-1)^jq^{-j^2}(j+1)!_{q^2}(j)!_{q^2}.$ For any $t\in \mathbb{Z}_{2p}$, assume that $$\Delta_{\mathbf{\overline{u}}_q(\mathfrak{sl_2})^*}\overline{(EK^t)}=\sum_{i,j,i',j'=0}^{p-1}
\sum_{l,l'\in\mathbb{Z}_{2p} } {^t\theta^{i',j',l'}_{i,j,l}}\overline{E^iF^jK^l}\ot \overline{E^{i'}F^{j'}K^{l'}}.$$
Then $$\begin{array}{rl}
{^t\theta^{i',j',l'}_{i,j,l}}=&\Delta_{\mathbf{\overline{u}}_q(\mathfrak{sl_2})^*}(\overline{EK^t})
({E^iF^jK^l}\ot {E^{i'}F^{j'}K^{l'}})\\
=&\overline{EK^t}(E^iF^jK^lE^{i'}F^{j'}K^{l'})\\
=&\overline{EK^t}(q^{2l(i'-j')}E^iF^jE^{i'}F^{j'}K^{l+l'}).\end{array}$$

{\bf Case 1}: $j=0$. In this case, $q^{2l(i'-j')}E^iF^jE^{i'}F^{j'}K^{l+l'}=
q^{2l(i'-j')}E^{i+i'}F^{j'}K^{l+l'}$.
Hence $$\begin{array}{rl}
{^t\theta^{i',j',l'}_{i,j,l}}
=&q^{2l(i'-j')}\delta_{i+i',1}\delta_{j',0}\delta_{l+l',t}\\
=&\left\{
\begin{array}{ll}
\vspace{0.2cm}
q^{2l}& {\rm if}\ i=0,\ i'=1 \ j'=0 \ {\rm and} \ l'=t-l\ in \ \mathbb{Z}_{2p} ,\\
1, & {\rm if} \ i=1,i'=0,j'=0 \ {\rm and} \ l'=t-l \ in \ \mathbb{Z}_{2p},\\
0, & otherwise.
\end{array}\right.
\end{array}$$
{\bf Case 2}: $1\<j\<i'$. In this case, by Lemma \ref{gs1}, we have
$$\begin{array}{rl}
q^{2l(i'-j')}E^iF^jE^{i'}F^{j'}K^{l+l'}
&=q^{2l(i'-j')}E^i(\sum_{s=0}^{j}E^{i'-s}F^{j-s}f^{i',j}_{s}(K))F^{j'}K^{l+l'}\\
&=\sum_{s=0}^{j}q^{2l(i'-j')}E^{i+i'-s}F^{j+j'-s}K^{l+l'}\phi^{i',j}_{s,j'}(K),\end{array}$$ where $\phi^{i',j}_{s,j'}(K)\in \Bbbk[K,K^{-1}]$ satisfying $f^{i',j}_{s}(K)F^{j'}=F^{j'}\phi^{i',j}_{s,j'}(K)$. Clearly, if $j'=0$ then $\phi^{i',j}_{s,0}(K)=f^{i',j}_s(K)$ for all $0\<s\<j.$ By $0\<s\<j,$ it is easy to see that if ${^t\theta^{i',j',l'}_{i,j,l}}\neq 0$, then $j'=0$ and $i+i'-j=1$. Since $1\<j\<i'<p$, $i+i'-j=1$ if $i=1$ and $i'=j$ or $i=0$ and $i'=j+1$. In the later case, $p>2$.

In case $j'=0$, $i=1$ and $i'=j$, we have
$$\begin{array}{rl}
{^t\theta^{j,0,l'}_{1,j,l}}
=&\overline{EK^t}(\sum_{s=0}^jq^{2lj}E^{1+j-s}F^{j-s}K^{l+l'}\phi^{j,j}_{s,0}(K))\\
=&\overline{EK^t}(q^{2lj}EK^{l+l'}f^{j,j}_j(K))\\
=&\overline{EK^t}(q^{2lj}\sum_{s=0}^jp^j_sEK^{l+l'+j-2s})\\
=&\sum_{s=0}^{j}q^{2lj}p^j_s\delta_{t,l+l'+j-2s}\\
=&\sum_{s=0}^{j}q^{2lj}p^j_s\delta_{l',t-l-j+2s}\\\end{array}$$
In case $j'=0$, $i=0$ and $i'=j+1$, we have
$$\begin{array}{rl}
{^t\theta^{j+1,0,l'}_{0,j,l}}
=&\overline{EK^t}(q^{2l(j+1)}EK^{l+l'}\phi^{j+1,j}_{j,0}(K))\\
=&\overline{EK^t}(q^{2l(j+1)}EK^{l+l'}f^{j+1,j}_j(K))\\
=&q^{2l(j+1)}\overline{EK^t}(\sum_{s=0}^j\varphi^j_sEK^{l+l'+j-2s})\\
=&q^{2l(j+1)}\sum_{s=0}^{j}\varphi^j_s\delta_{t,l+l'+j-2s}\\
=&\sum_{s=0}^{j}q^{2l(j+1)}\varphi^j_s\delta_{l',t-l-j+2s}\\\end{array}$$

{\bf Case 3}: $j>i'$,. In this case, by Lemma \ref{gs2}, we have
$$\begin{array}{rl}
q^{2l(i'-j')}E^iF^jE^{i'}F^{j'}K^{l+l'}
&=q^{2l(i'-j')}E^i(\sum_{s=0}^{i'}E^{i'-s}F^{j-s}g^{i',j}_s(K))F^{j'}K^{l+l'}\\
&=\sum_{s=0}^{i'}q^{2l(i'-j')}E^{i+i'-s}F^{j+j'-s}K^{l+l'}\psi^{i',j}_{s,j'}(K)\\
\end{array},$$
 where $\psi^{i',j}_{s,j'}(K)\in \Bbbk[K,K^{-1}]$ satisfying $g^{i',j}_s(K)F^{j'}=F^{j'}\psi^{i',j}_{s,j'}(K)$. By $j>i'$, $j+j'-s\>1$ for all $0\<s\<i'$. Hence ${^t\theta^{i',j',l'}_{i,j,l}}=0$ if $j>i'$.

 By the discussion above, we have
 $$\begin{array}{rl}
 &\Delta_{\mathbf{\overline{u}}_q(\mathfrak{sl_2})^*}(\a)\\
 =&\sum_{t\in \mathbb{Z}_{2p}}\Delta_{\mathbf{\overline{u}}_q(\mathfrak{sl_2})^*}(\overline{EK^t})\\
 =&\sum_{t\in  \mathbb{Z}_{2p}}(\sum_{l\in  \mathbb{Z}_{2p}}q^{2l}\overline{K^l}\ot \overline{EK^{l-t}}+\sum_{l\in  \mathbb{Z}_{2p}}\overline{EK^l}\ot \overline{K^{l-t}})\\
 &+\sum_{j=1}^{p-1}\sum_{l,l' \in \mathbb{Z}_{2p}}\sum_{s=0}^jq^{2lj}p^j_s\delta_{l',t-l-j+2s}\overline{EF^jK^l}\ot \overline{E^jK^{l'}}\\
 &+\sum_{1\<j\<p-2}\sum_{l,l'\in \mathbb{Z}_{2p}}\sum_{s=0}^jq^{2l(j+1)}\varphi^j_s\delta_{l',t-l-j+2s}\overline{F^jK^l}\ot \overline{E^{j+1}K^{l'}}\\
 =&(\sum_{l\in \mathbb{Z}_{2p}}q^{2l}\overline{K^l})\ot (\sum_{l\in \mathbb{Z}_{2p}}\overline{EK^t})
 +(\sum_{l\in \mathbb{Z}_{2p}}\overline{EK^l})\ot (\sum_{l\in \mathbb{Z}_{2p}}\overline{K^t})\\
 &+\sum_{j=1}^{p-1}\sum_{t,l\in \mathbb{Z}_{2p}}\sum_{s=0}^jq^{2lj}p^{j}_s\overline{EF^jK^l}\ot\overline{E^jK^{t-l-j+2s}}\\
 &+\sum_{j=1}^{p-2}\sum_{t,l\in \mathbb{Z}_{2p}}\sum_{s=0}^jq^{2l(j+1)}\varphi^{j}_s\overline{F^jK^l}\ot\overline{E^{j+1}K^{t-l-j+2s}}\\
 =&\g^2\ot \a+\a\ot \varepsilon +\sum_{j=1}^{p-1}\sum_{l\in \mathbb{Z}_{2p}}\sum_{s=0}^jq^{2lj}p^j_s\overline{EF^jK^l}\ot(\sum_{t\in \mathbb{Z}_{2p}}\overline{E^jK^t})\\
 &+\sum_{1\<j\<p-2}\sum_{l\in \mathbb{Z}_{2p}}q^{2l(j+1)}(\sum_{s=0}^j\varphi^j_s)\overline{F^jK^l}\ot (\sum_{t\in \mathbb{Z}_{2p}}\overline{E^{j+1}K^t})\\
 =&\g^2\ot \a+\a\ot \varepsilon+\sum_{1\<j\<p-2}(-1)^jq^{-j^2}\b^j\g^{2(j+1)}\ot \a^{j+1}\\
 =&\a \ot \varepsilon +\sum_{j=0}^{p-2}(-1)^{j}q^{-j^2}\b^j\g^{2(j+1)}\ot \a^{j+1}.\\
 \end{array}$$
The identities \textup{(b)} and \textup{(c)} are proved similarly.
\end{proof}

\begin{proposition}\label{2.10}\label{2.11}\label{2.12}
In $\uq^{*\cop}$, the following identities hold:
\begin{enumerate}
\item[(a)]
$S(\a)=\sum_{j=0}^{p-2}(-1)^{j+1}q^{j^2+4j}\a^{j+1}\b^j \g^{-2}.$
\item[(b)]
$S(\b)=\sum_{i=0}^{1}-q^{i}\a^{i}\b^{i+1} \g^{2}.$
\item[(c)]
$S(\g)=\sum_{j=0}^{p-1}(-1)^jq^{j^2+2j}\a^{j}\b^{j} \g^{-1}.$
\end{enumerate}
\end{proposition}
\begin{proof}
\textup{(a)}
For any fixed $t\in \mathbb{Z}_{2p}$, assume $S(\overline{EK^t})=\sum_{i,j=0}^{p-1}\sum_{l\in \mathbb{Z}_{2p}}{^t\theta_{i,j,l}}\overline{E^iF^jK^l}.$  Then $$\begin{array}{rl}
{^t\theta_{i,j,l}}=S(\overline{EK^t})(E^iF^jK^l)=(S^{-1})^*\overline{EK^t}(E^iF^jK^l)=
\overline{EK^t}(S^{-1}(E^iF^jK^l))\end{array}$$ One can easy check that
$S^{-1}(E^iF^jK^l)=(-1)^{i+j}q^{(j-i)(1-j+i+2l)}F^jE^iK^{j-i-l}.$

Case 1: If $i=0$, then one can check $\overline{EK^t}(S^{-1}(E^iF^jK^l))=0$. Thus,${^t\theta_{0,j,l}}=0$.

Case 2: If $i\>1$ and $i=j$, then by Lemma \ref{gs1}, we have $$\overline{EK^t}(S^{-1}(E^iF^jK^l))=\overline{EK^t}(\sum_{y=0}^jE^{i-y}F^{j-y}f^{i,j}_{y}(K)K^{-l}).$$ It is easy to see $\overline{EK^t}(S^{-1}(E^iF^jK^l))=0$. Thus ${^t\theta_{i,j,l}}=0$ if $i=j\>1$.

Case 3: If $i\>1$ and $i<j$. Then by Lemma \ref{gs2}, we have $$\begin{array}{rl}
{^t\theta_{i,j,l}}=&\overline{EK^t}(S^{-1}(E^iF^jK^l))\\
=&(-1)^{i+j}q^{(j-i)(1-j+i+2l)}\overline{EK^t}(\sum_{y=0}^i
E^{i-y}F^{j-y}g^{i,j}_{y}(K)K^{j-i-l}).\end{array}$$ It is easy to see $\overline{EK^t}(S^{-1}(E^iF^jK^l))=0$. Thus ${^t\theta_{i,j,l}}=0$ if $i\>1$ and $i<j$.

Case 4: If $i\>1$ and $i>j$.

 In case $j=0$. We have $$\begin{array}{rl}
{^t\theta_{i,0,l}}
=&\overline{EK^t}((-1)^iq^{-i(1+i+2l)}E^iK^{-i-l})\\
=&(-1)^iq^{-i(1+i+2l)}\overline{EK^t}(E^iK^{-i-l})\\
=&(-1)^iq^{-i(1+i+2l)}\delta_{i,1}\delta_{t,-i-l}.\\\end{array}$$

In case $j\>1$. Then $$\begin{array}{rl}
{^t\theta_{i,j,l}}=&\overline{EK^t}(S^{-1}(E^iF^jK^l))\\
=&(-1)^{i+j}q^{(j-i)(1-j+i+2l)}
\overline{EK^t}(\sum_{y=0}^jE^{i-y}F^{j-y}f^{i,j}_{y}(K)K^{j-i-l}))
\end{array}$$
If ${^t\theta_{i,j,l}}\neq0$, then one gets $i=j+1$. In this case, we have
$$\begin{array}{rl}
{^t\theta_{j+1,j,l}}=&\overline{EK^t}(S^{-1}(E^iF^jK^l))\\
=&-q^{-2(l+1)}\overline{EK^t}(Ef^{j+1,j}_j(K)K^{-1-l})\\
=&-q^{-2(l+1)}(\sum_{s=0}^j\varphi^j_s\delta_{t,j-l-1-2s}),
\end{array}$$
where $\varphi^j_y$ are defined in Proposition \ref{2.7} for $0\<y\<j$.

 By the discussion above, we have
 $$\begin{array}{rl}
 S(\a)=&(\sum_{t\in \mathbb{Z}_{2p}}-q^{2t}\overline{EK^{-1-t}})\\
 &+\sum_{t\in \mathbb{Z}_{2p}}\sum_{j=1}^{p-2}\sum_{y=0}^j-q^{-2(j-t-2y)}
 \varphi^j_y\overline{E^{j+1}F^jK^{j-t-1-2y}}.
 \end{array}$$
 By a straightforward computation, one gets $(\sum_{t\in \mathbb{Z}_{2p}}-q^{2t}\overline{EK^{-1-t}})=-\a\g^{-2}$ and
$$\begin{array}{rl}
&\sum_{t\in \mathbb{Z}_{2p}}\sum_{j=1}^{p-2}\sum_{y=0}^j-q^{-2(j-t-2y)}\varphi^j_y
\overline{E^{j+1}F^jK^{j-t-1-2y}}\\
=&\sum_{j=1}^{p-2}\sum_{t\in\mathbb{Z}_{2p}}-q^{2(t-j)}(\sum_{y=0}^j\varphi^j_y)\overline{E^{j+1}F^jK^{j-t-1}},\\
=&\sum_{j=1}^{p-2}(-1)^{j+1}q^{j^2+4j}\a^{j+1}\b^j\g^{-2} \end{array}$$

Consequently, we have
$$S(\a)=\sum_{j=0}^{p-2}(-1)^{j+1}q^{j^2+4j}\a^{j+1}\b^j \g^{-2}$$
The identities \textup{(b)} and \textup{(c)} are proved similarly.
\end{proof}

\begin{lemma}\label{2.13}
The left $\mathbf{\overline{u}}_q(\mathfrak{sl_2})$-module structure of $\mathbf{\overline{u}}_q(\mathfrak{sl_2})^*$ is determined by
$$\begin{array}{rcl}
K\rightharpoonup \a=\a, &K\rightharpoonup \b=\b, &K\rightharpoonup \g=q\g,\\
E\rightharpoonup \a=\g^2, &E\rightharpoonup \b=q^{-1}\b^2\g^2, &E\rightharpoonup \g=-\b\g^3,\\
F\rightharpoonup \a=0, &F\rightharpoonup \b=\g^{-2}, &F\rightharpoonup \g=0.\\
\end{array}$$
\end{lemma}
\begin{proof}
It follows from a straightforward verification.
\end{proof}

\begin{lemma}\label{2.14}
The right $\mathbf{\overline{u}}_q(\mathfrak{sl_2})$-module structure of $\mathbf{\overline{u}}_q(\mathfrak{sl_2})^*$ is determined by
$$\begin{array}{rcl}
\a\leftharpoonup K=q^2\a, &\b \leftharpoonup K=q^{-2}\b, &\g\leftharpoonup K=q\g,\\
\a\leftharpoonup E=\varepsilon, &\b\leftharpoonup E=0, &\g\leftharpoonup E=0,\\
\a\leftharpoonup F=-q^{-1}\a^2, &\b\leftharpoonup F=\varepsilon+q(2)_{q^2}\a\b, &\g\leftharpoonup F=-q^{-1}\a\g.\\
\end{array}$$
\end{lemma}
\begin{proof}
It follows from a straightforward verification.
\end{proof}
Since $\uq$ is generated by $K,E,F$, the preceding two lemmas determine the left and right $\uq$-module structures on $\uq^*$.
\begin{theorem}\label{2.16}
The Drinfeld double $D(\mathbf{\overline{u}}_q(\mathfrak{sl_2}))$ is generated by $E, F, K$, $\a,\b,\g$ subject to the relations defining $\mathbf{\overline{u}}_q(\mathfrak{sl_2})$ and $\mathbf{\overline{u}}_q(\mathfrak{sl_2})^*$ and the following relations:
\begin{enumerate}
\item[(a)]
$E \a=q^{-2}\a E+q^{-2}\g^2K-q^{-2}$
\item[(b)]
$E\b=q^2\b E+q\b^2\g^2K$
\item[(c)]
$E\g=q^{-1}\g E-q^{-1}\b\g^3 K$

\item[(d)]
$[F,\a]=q^{3}\a^2K^{-1}$
\item[(e)]
$[F,\b]=\g^{-2}-K^{-1}-q^2[2]\a\b K^{-1}$
\item[(f)]
$F\g=q\a\g K^{-1}+q^{-1}\g F$

\item[(g)]
$K\a=q^{-2}\a K$
\item[(h)]
$K \b=q^2\b K$
\item[(i)]
$K\g=\g K$
\end{enumerate}
\end{theorem}
\begin{proof} We only prove \textup{(a)} since the proof of others are similar.
Put $U=\uq$. We compute \textup{(a)} directly from the cross
multiplication formula
\[
 hf=\sum_{(h)}\bigl(h_1\rightharpoonup f
                \leftharpoonup S^{-1}(h_3)\bigr)h_2,
 \qquad h\in U,\quad f\in U^*.
\]
Taking $h=E$, we use\[
 \Delta^{(2)}(E)=E\otimes K\otimes K
       +1\otimes E\otimes K+1\otimes1\otimes E,
\]
and
\[
 S^{-1}(K)=K^{-1},\qquad S^{-1}(E)=-K^{-1}E.
\]
From Lemma~\ref{2.14},
\[
 \alpha\leftharpoonup K^{-1}=q^{-2}\alpha,\qquad
 \beta\leftharpoonup K^{-1}=q^2\beta,\qquad
 \gamma\leftharpoonup K^{-1}=q^{-1}\gamma.
\]
We also use
\[
 (f\leftharpoonup a)\leftharpoonup b=f\leftharpoonup(ab),
 \qquad
 (fg)\leftharpoonup K^{-1}
 =(f\leftharpoonup K^{-1})(g\leftharpoonup K^{-1}),
\]
where the second identity holds because $K^{-1}$ is group-like.

(a)
Substituting $h=E$ and $f=\alpha$ directly into the cross multiplication,
and then using Lemmas~\ref{2.13}--\ref{2.14}, gives
\begin{align*}
 E\alpha
 &=\sum_{(E)}\bigl(E_1\rightharpoonup\alpha
                      \leftharpoonup S^{-1}(E_3)\bigr)E_2\\
 &=\bigl((E\rightharpoonup\alpha)\leftharpoonup K^{-1}\bigr)K
   +(\alpha\leftharpoonup K^{-1})E
   -\alpha\leftharpoonup(K^{-1}E)\\
 &=(\gamma^2\leftharpoonup K^{-1})K
   +q^{-2}\alpha E
   -(\alpha\leftharpoonup K^{-1})\leftharpoonup E\\
 &=(q^{-1}\gamma)^2K+q^{-2}\alpha E
   -q^{-2}(\alpha\leftharpoonup E)\\
 &=q^{-2}\gamma^2K+q^{-2}\alpha E-q^{-2}\varepsilon\\
 &=q^{-2}\alpha E+q^{-2}\gamma^2K-q^{-2}.
\end{align*}
In the last line, the unit $\varepsilon$ of $U^*$ is identified with
$1$ in $D(U)$. This proves \textup{(a)}.

Lastly, these relations give a presentation, not merely identities in the double. The canonical map from the algebra with the stated generators and relations onto $D(U)$ is surjective. The cross relations move each occurrence of a generator of $U$ past a word from $U^{*\cop}$, so this presented algebra is spanned by the products of the $2p^3$ basis elements of $U^{*\cop}$ with the $2p^3$ basis elements of $U$. Its dimension is at most $4p^6=\dim D(U)$; surjectivity forces equality and hence the map is an isomorphism.
\end{proof}

In the following of this paper, denote $D(\uq)$ by $H$.
\section{Irreducible representations of $H$}\selabel{4}
In this section, $p\ge2$ is arbitrary and $q$ has order $2p$. Put
$d=q-q^{-1}$.
We use the generators and relations of Theorem~\ref{2.16}. In particular, the algebra underlying $\uq^{*\cop}$ is the algebra described in Lemmas 2.2--2.3. For $t\in\Z$, let $\res t$ be its least nonnegative residue modulo $p$. Representatives of $\Z_{2p}$ are taken in $\{0,1,\ldots,2p-1\}$.

For $(i,j)\in\Z_{2p}^2$, define
\begin{align}
 r(i,j)&=\begin{cases}p,&i\text{ odd},\\1+\res{j+i/2},&i\text{ even},\end{cases}\label{eq:r}\\
 s(i,j)&=\begin{cases}p,&i\text{ odd},\\1+\res{i/2-j},&i\text{ even}.\end{cases}\label{eq:s}
\end{align}
The formulas are well defined: replacing an even representative $i$ by $i+2p$ changes $i/2$ by $p$, and parity is well defined modulo $2p$. We first develop the change of generators and the highest weight construction, and then prove the classification at the end of this section.

\subsection{Modified generators and a PBW basis}
Set
\begin{equation}
 A=\alpha K^{-1},\quad T=\beta\gamma^2K,\quad c_0=\frac{q^2}{d},
 \qquad X=E+\frac Td,\quad Y=F-c_0A.
 \label{eq:change}
\end{equation}
Here $d\ne0$ because $q$ has order $2p\ge4$.
\begin{proposition}\label{prop:generators}
The algebra $H$ is generated by $X,Y,K,\alpha,\beta,\gamma$.
\end{proposition}
\begin{proof}
Let $H'$ be the unital subalgebra generated by the six indicated elements. Since $K^{-1}=K^{2p-1}\in H$, the definitions show that $X,Y\in H$, and hence $H'\subseteq H$. Conversely, $K^{-1}=K^{2p-1}\in H'$ because $K\in H'$. Consequently $\beta\gamma^2K$ and $\alpha K^{-1}$ belong to $H'$, and
\begin{equation}
 E=X-d^{-1}\beta\gamma^2K,\qquad
 F=Y+q^2d^{-1}\alpha K^{-1}
 \label{eq:recover}
\end{equation}
belong to $H'$. The original generators $E,F,K,\alpha,\beta,\gamma$ are therefore all contained in $H'$, proving $H\subseteq H'$. This completes the proof.
\end{proof}

\begin{lemma}\label{lem:modified}
The modified generators satisfy
\begin{align}
 Y\alpha&=\alpha Y, & XY&=YX, & X\beta&=q^2\beta X,\label{eq:commute}\\
 X\alpha&=q^{-2}\alpha X+q^{-2}(\gamma^2K-1),\label{eq:xa}\\
 \beta Y&=Y\beta+K^{-1}-\gamma^{-2},\label{eq:by}\\
 KX&=q^2XK,&KY&=q^{-2}YK,\label{eq:Kweight}\\
 \gamma X&=qX\gamma,&\gamma Y&=qY\gamma.\label{eq:Gweight}
\end{align}
Moreover, $X^p=Y^p=0$.
\end{lemma}
\begin{proof}
 The original relations give $K^{-1}\alpha=q^2\alpha K^{-1}$ and therefore
\[
 [Y,\alpha]=\bigl(q^3-c_0(q^2-1)\bigr)\alpha^2K^{-1}=0.
\]
Rewriting the relation for $F\gamma$ yields
 $\gamma F=qF\gamma-q^2\alpha\gamma K^{-1}.$
Since $\gamma A=q^{-1}A\gamma$, it follows that
\[
 \gamma Y-qY\gamma=
 \bigl(-q^2+c_0(q-q^{-1})\bigr)A\gamma=0.
\]
The relation $KY=q^{-2}YK$ follows directly from the weights of $F$ and $A$. Also,
\begin{align*}
 [\beta,Y]
 &=K^{-1}-\gamma^{-2}
 +\bigl(q^2[2]-c_0(q^2-q^{-2})\bigr)\alpha\beta K^{-1}\\
 &=K^{-1}-\gamma^{-2},
\end{align*}
where $q^2-q^{-2}=d[2]$.
For $X$, observe that
\[
 T\beta=\beta T,\quad T\alpha=q^{-2}\alpha T,
 \quad KT=q^2TK,\quad\gamma T=q^{-1}T\gamma.
\]
It follows that
\[
 X\beta-q^2\beta X
 =\left(q+\frac{1-q^2}{d}\right)\beta^2\gamma^2K=0.
\]
The relation for $X\alpha$ follows from Theorem~\ref{2.16}(a) and $T\alpha=q^{-2}\alpha T$. The relation for $KX$ is immediate, while
\[
 X\gamma-q^{-1}\gamma X
 =\left(-q^{-1}+\frac{1-q^{-2}}{d}\right)\beta\gamma^3K=0.
\]
To prove $XY=YX$, use $K^{-1}E=q^{-2}EK^{-1}$ to obtain
\[
 [E,A]=q^{-2}(\gamma^2-K^{-1}),\qquad
 [E,Y]=\frac{K-\gamma^2}{d}.
\]
The already established weight relations imply that $\gamma^2K$ commutes with $Y$. Hence
\[
 [T,Y]=[\beta,Y]\gamma^2K=\gamma^2-K,
\]
and $[X,Y]=[E,Y]+d^{-1}[T,Y]=0$.

It remains to prove nilpotence. The above relations imply
\[
 YA=q^{-2}AY,\qquad XT=q^{-2}TX.
\]
Since $q^{-2}$ is a primitive $p$-th root of unity, all intermediate Gaussian coefficients in the $p$-th power vanish. Thus the quantum binomial formula gives
\[
 F^p=Y^p+c_0^pA^p,
 \qquad E^p=X^p+(-d^{-1})^pT^p.
\]
The elements $A^p$ and $T^p$ are scalar multiples of $\alpha^pK^{-p}$ and $\beta^p\gamma^{2p}K^p$, respectively, and are zero. Since $E^p=F^p=0$, we obtain $X^p=Y^p=0$.
\end{proof}

\begin{lemma}\label{lem:pbw}
The following set is a basis of $H$:
\begin{equation}
 \{\alpha^aY^bK^k\gamma^\ell X^c\beta^e:
 0\le a,b,c,e<p,\ k,\ell\in\Z_{2p}\}.
 \label{eq:pbw}
\end{equation}
The subalgebra $B=\langle K,\gamma,X,\beta\rangle$ has basis
\[
 \{K^k\gamma^\ell X^c\beta^e:
 k,\ell\in\Z_{2p},\ 0\le c,e<p\}.
\]
In particular, $H$ is a free right $B$-module with basis $\{\alpha^aY^b:0\le a,b<p\}$.
\end{lemma}
\begin{proof}
By Proposition~\ref{prop:generators}, $H$ is generated by
$\alpha,Y,K,\gamma,X,\beta$. By the defining relations of $\uq$, Lemmas~\ref{2.2}
and~\ref{lem:modified}, and Theorem~\ref{2.16}, every word in these generators can be
rewritten as a linear combination of the ordered monomials in
\eqref{eq:pbw}. Hence the family
in \eqref{eq:pbw} spans $H$. It has $p^4(2p)^2=4p^6$ members,
whereas
\[
 \dim H=(\dim\uq)^2=(2p^3)^2=4p^6.
\]
Comparison of dimensions shows that \eqref{eq:pbw} is a basis.

The same relations show that $B$ is spanned by
$K^k\gamma^\ell X^c\beta^e$ with the stated ranges. These monomials
form a subfamily of \eqref{eq:pbw}, so they are linearly independent
and therefore form a basis of $B$. Grouping the basis
\eqref{eq:pbw} according to the left factor $\alpha^aY^b$ gives
\[
 H=\bigoplus_{0\le a,b<p}\alpha^aY^bB.
\]
Thus $H$ is a free right $B$-module with the asserted basis.
\end{proof}
 In fact one can check that the presentation of $H$ are given by  Lemmas \ref{2.2}, \ref{lem:modified}
 and \[\quad K^{2p}=1,\quad K\gamma=\gamma K,\quad
 K\alpha=q^{-2}\alpha K,\quad K\beta=q^2\beta K.
\]
\subsection{Highest weight induced modules}
For an $H$-module $M$, put
\[
 M^X=\{m\in M:Xm=0\},\qquad M^\beta=\{m\in M:\beta m=0\}.
\]
The space $M^X\cap M^\beta$ equals $M^E\cap M^\beta$, since $X-E=d^{-1}\beta\gamma^2K$ annihilates $M^\beta$.

\begin{lemma}\label{lem:highest}
Every nonzero $H$-module contains a nonzero vector $v$ such that
\begin{equation}
 Xv=\beta v=0,\qquad Kv=q^iv,\qquad\gamma v=q^jv
 \label{eq:highest}
\end{equation}
for some $(i,j)\in\Z_{2p}^2$.
\end{lemma}
\begin{proof}
Since $X^p=0$, choose $0\ne m\in M^X$. Let $t\ge1$ be minimal with $\beta^tm=0$. Then $v'=\beta^{t-1}m\ne0$ is killed by $\beta$, and $Xv'=q^{2(t-1)}\beta^{t-1}Xm=0$. Thus $N=M^X\cap M^\beta$ is nonzero. The weight relations show that $N$ is stable under the commuting operators $K,\gamma$. Hence there is a common eigenvector $v$ such that  \eqref{eq:highest}.
\end{proof}

For $(i,j)\in\Z_{2p}^2$, let $\kk_{i,j}$ be the one-dimensional $B$-module defined by
\[
 X\cdot1=\beta\cdot1=0,\qquad K\cdot1=q^i,\qquad\gamma\cdot1=q^j.
\]
This is well defined. Indeed, the ordered basis of $B$ shows that sending $X,\beta$ to zero gives a quotient isomorphic to the group algebra $\kk[\Z_{2p}^2]$. Define
\begin{equation}
 Z(i,j)=H\otimes_B\kk_{i,j},\qquad v=1\otimes1.
 \label{eq:induced}
\end{equation}
By Lemma~\ref{lem:pbw}, $Z(i,j)$ has basis
\begin{equation}
 z_{a,b}=\alpha^aY^bv,\qquad 0\le a,b<p,
 \label{eq:zbase}
\end{equation}
and hence has dimension $p^2$.

\begin{lemma}\label{lem:action}
On the basis \eqref{eq:zbase}, the generator actions are
\begin{align}
 Kz_{a,b}&=q^{i-2a-2b}z_{a,b},&
 \gamma z_{a,b}&=q^{j-a+b}z_{a,b},\label{eq:zweights}\\
 \alpha z_{a,b}&=z_{a+1,b},&Yz_{a,b}&=z_{a,b+1},\label{eq:zraise}\\
 Xz_{a,b}&=C_a z_{a-1,b},&\beta z_{a,b}&=D_{a,b}z_{a,b-1},\label{eq:zlower}
\end{align}
where vectors outside $0\le a,b<p$ are zero, $C_0=D_{a,0}=0$, and
\begin{align}
 C_a&=q^{-2}(a)_{q^{-2}}\bigl(q^{2j+i-2(a-1)}-1\bigr),\label{eq:Ca}\\
 D_{a,b}&=q^{2a}(b)_{q^2}\bigl(q^{-i}-q^{-2j-2(b-1)}\bigr).\label{eq:Dab}
\end{align}
\end{lemma}
\begin{proof}
The diagonal actions and the actions of $\alpha,Y$ follow directly from Thorem \ref{2.16} and Lemma~\ref{lem:modified}. Put $D=\gamma^2K$. Then $D\alpha=q^{-4}\alpha D$ and $DY=YD$. Repeated use of \eqref{eq:xa} gives, for $a\ge1$,
\begin{align*}
 X\alpha^a
 &=q^{-2a}\alpha^aX
 +q^{-2}\alpha^{a-1}\sum_{t=0}^{a-1}q^{-2t}
       \bigl(q^{-4(a-1-t)}D-1\bigr)\\
 &=q^{-2a}\alpha^aX
 +q^{-2}(a)_{q^{-2}}\alpha^{a-1}
       \bigl(q^{-2(a-1)}D-1\bigr).
\end{align*}
For the second equality, use
\[
 \sum_{t=0}^{a-1}q^{-2t}q^{-4(a-1-t)}
   =q^{-2(a-1)}(a)_{q^{-2}}.
\]
Applying this identity to $Y^bv$, and using $XY^bv=0$ and $DY^bv=q^{2j+i}Y^bv$, proves \eqref{eq:Ca}.

Similarly, \eqref{eq:by} and
$K^{-1}Y=q^2YK^{-1}$, $\gamma^{-2}Y=q^{-2}Y\gamma^{-2}$ imply
\[
 \beta Y^b=Y^b\beta+
 Y^{b-1}\bigl((b)_{q^2}K^{-1}-(b)_{q^{-2}}\gamma^{-2}\bigr).
\]
Now $\beta\alpha^a=q^{2a}\alpha^a\beta$, and
$(b)_{q^{-2}}=q^{-2(b-1)}(b)_{q^2}$. Applying the identity to $v$ yields \eqref{eq:Dab}.
\end{proof}

\begin{lemma}\label{lem:firstzero}
Let $r=r(i,j)$ and $s=s(i,j)$ be as in \eqref{eq:r}--\eqref{eq:s}. Then $C_a\ne0$ for $1\le a<r$; if $r<p$, then $C_r=0$. Also, $D_{a,b}\ne0$ for $0\le a<p$ and $1\le b<s$; if $s<p$, then $D_{a,s}=0$ for all $a$.
\end{lemma}
\begin{proof}
Recall from \eqref{eq:Ca}--\eqref{eq:Dab} that
\[
 C_a=q^{-2}(a)_{q^{-2}}
       \bigl(q^{2j+i-2(a-1)}-1\bigr),\qquad
 D_{a,b}=q^{2a}(b)_{q^2}
       \bigl(q^{-i}-q^{-2j-2(b-1)}\bigr).
\]
Since $q$ has order $2p$, both $q^2$ and $q^{-2}$ have order $p$.
For either choice $\xi=q^2$ or $q^{-2}$, and for $1\le n<p$, we have
\[
 (n)_\xi=\frac{1-\xi^n}{1-\xi}\ne0,
\]
because $\xi\ne1$ and $\xi^n\ne1$. Moreover, every integral power of
$q$ is nonzero. Thus, in the ranges $1\le a<p$ and $1\le b<p$,
respectively, the coefficients can vanish only through the last
parenthesized factors in the displayed formulas.

We first determine the zeros of $C_a$. The order of $q$ gives
\begin{align*}
 C_a=0
 &\quad\Longleftrightarrow\quad q^{2j+i-2(a-1)}=1\\
 &\quad\Longleftrightarrow\quad
       2j+i-2(a-1)\equiv0\pmod{2p}.
\end{align*}
If $i$ is odd, the integer $2j+i-2(a-1)$ is odd and cannot be
divisible by $2p$. Hence $C_a\ne0$ for every $1\le a<p$.
In this case $r=p$ by \eqref{eq:r}, which proves the assertion about
$C_a$.

Suppose now that $i$ is even, and write $i=2t$ with $t\in\mathbb Z$.
Dividing the congruence by $2$ changes its modulus from $2p$ to $p$;
therefore
\[
 C_a=0\quad\Longleftrightarrow\quad
 a-1\equiv j+t\pmod p.
\]
Put $\rho=\res{j+t}\in\{0,1,\ldots,p-1\}$. Then
$r=1+\rho$ by \eqref{eq:r}. If $1\le a<r$, then
$0\le a-1<\rho<p$, so $a-1$ cannot be congruent to $\rho$ modulo $p$.
Consequently $C_a\ne0$ throughout this range. If $r<p$, then
$1\le r<p$ and $r-1=\rho$, so the same equivalence gives $C_r=0$.

We next determine the zeros of $D_{a,b}$. For every fixed
$0\le a<p$, the factor $q^{2a}$ is nonzero, so the zero condition
is independent of $a$. For $1\le b<p$, we obtain
\begin{align*}
 D_{a,b}=0
 &\quad\Longleftrightarrow\quad
       q^{-i}=q^{-2j-2(b-1)}\\
 &\quad\Longleftrightarrow\quad
       2(b-1)\equiv i-2j\pmod{2p}.
\end{align*}
If $i$ is odd, this congruence has no solution: its left-hand side
is even, whereas its right-hand side is odd modulo the even
modulus $2p$. Thus $D_{a,b}\ne0$ for all $0\le a<p$ and
$1\le b<p$. Since $s=p$ in this case, the required assertion follows.

If $i=2t$ is even, the congruence is equivalent to
\[
 b-1\equiv t-j\pmod p.
\]
Put $\sigma=\res{t-j}\in\{0,1,\ldots,p-1\}$, so that
$s=1+\sigma$ by \eqref{eq:s}. If $1\le b<s$, then
$0\le b-1<\sigma<p$, which rules out this congruence.
It follows that $D_{a,b}\ne0$ for every $0\le a<p$ in the asserted
range. If $s<p$, then $s-1=\sigma$, so $D_{a,s}=0$ for every $a$.
This proves all the claims.
\end{proof}

\subsection{Simple quotients and completeness}
For $r=r(i,j)$ and $s=s(i,j)$, set
\begin{equation}
 W(i,j)=\Span\{z_{a,b}:a\ge r\text{ or }b\ge s\}\subseteq Z(i,j).
 \label{eq:W}
\end{equation}
Only indices from the basis \eqref{eq:zbase} are included. In particular, $r=p$ or $s=p$ contributes no excluded vectors in that direction.

\begin{proposition}\label{prop:simple}
The space $W(i,j)$ is an $H$-submodule. The quotient
\[
 L(i,j)=Z(i,j)/W(i,j)
\]
is simple, with basis
\[
 \{v_{a,b}=z_{a,b}+W(i,j):0\le a<r,\ 0\le b<s\}.
\]
\end{proposition}
\begin{proof}
Write $Z=Z(i,j)$, $W=W(i,j)$, $r=r(i,j)$ and $s=s(i,j)$.
We first check that $W$ is stable under each of the generators in
Proposition~\ref{prop:generators}. Throughout, a basis vector with an
index outside $\{0,\ldots,p-1\}$ is interpreted as zero.

\smallskip\noindent\emph{Stability of $W$.}
Let $z_{a,b}$ be a basis vector belonging to $W$, so that $a\ge r$
or $b\ge s$. By \eqref{eq:zweights}, both $Kz_{a,b}$ and
$\gamma z_{a,b}$ are scalar multiples of $z_{a,b}$ and belong to $W$.
By \eqref{eq:zraise}, the operators $\alpha$ and $Y$ only increase
one of the indices. Thus any nonzero image still satisfies
$a\ge r$ or $b\ge s$, and $\alpha W,YW\subseteq W$.

Next, $Xz_{a,b}=C_a z_{a-1,b}$. If $b\ge s$, its nonzero image
still has second index at least $s$, so lies in $W$. Otherwise
$b<s$, and membership of $z_{a,b}$ in $W$ forces $a\ge r$.
If $a>r$, then $a-1\ge r$, so again the image lies in $W$.
The remaining possibility is $a=r$. Such a basis vector can occur
only if $r<p$, since $a<p$; Lemma~\ref{lem:firstzero} then gives
$C_r=0$. Hence $Xz_{r,b}=0$. This proves $XW\subseteq W$.

Similarly, $\beta z_{a,b}=D_{a,b}z_{a,b-1}$. If $a\ge r$,
every nonzero image remains in $W$. If $a<r$, then $b\ge s$.
For $b>s$ we have $b-1\ge s$, whereas for $b=s$ the basis vector
can occur only when $s<p$, in which case
$D_{a,s}=0$ by Lemma~\ref{lem:firstzero}. Thus $\beta W\subseteq W$.
These arguments also cover $r=p$ or $s=p$, when the corresponding
boundary vectors do not occur. Since $X,Y,K,\alpha,\beta,\gamma$
generate $H$, stability under these six operators proves that $W$
is an $H$-submodule.

\smallskip\noindent\emph{The quotient basis.}
The basis \eqref{eq:zbase} is the disjoint union of the vectors
spanning $W$ and the vectors with $0\le a<r$, $0\le b<s$.
Consequently, as vector spaces,
\[
 Z=W\oplus\Span\{z_{a,b}:0\le a<r,\ 0\le b<s\}.
\]
The images $v_{a,b}=z_{a,b}+W$ of the latter vectors therefore form
a basis of $L=Z/W$. In particular, $v_{0,0}\ne0$ since $r,s\ge1$.
Hence $W$ is proper and $L$ is nonzero. The generator actions descend
to this quotient; in particular,
\[
 Xv_{a,b}=C_av_{a-1,b},\qquad
 \beta v_{a,b}=D_{a,b}v_{a,b-1},
\]
with an image outside the rectangle interpreted as zero.

\smallskip\noindent\emph{ Simplicity of $Z/W$.}
We prove directly that $W$ is a maximal proper $H$-submodule of $Z$.
Let $N$ be an $H$-submodule satisfying
\[
 W\subsetneq N\subseteq Z.
\]
Choose $w\in N\setminus W$. Using the direct-sum decomposition
above, write uniquely
\[
 w=w_W+\sum_{a=0}^{r-1}\sum_{b=0}^{s-1}\lambda_{a,b}z_{a,b},
 \qquad w_W\in W.
\]
Since $w\notin W$, at least one coefficient $\lambda_{a,b}$ is
nonzero. As $W\subseteq N$, we may subtract $w_W$ and obtain
\[
 w_0=\sum_{a=0}^{r-1}\sum_{b=0}^{s-1}
                  \lambda_{a,b}z_{a,b}\in N\setminus\{0\}.
\]
Thus the following calculation takes place in the induced module
$Z$, before any simplicity assertion about its quotient.

Choose
 $a_0=\max\{a:\lambda_{a,b}\ne0\text{ for some }b\},
 b_0=\max\{b:\lambda_{a_0,b}\ne0\}.$
The formula $Xz_{a,b}=C_az_{a-1,b}$ shows that $X^{a_0}$ kills
every term with $a<a_0$: after the first index reaches zero,
the next application gives zero because $C_0=0$. For $a=a_0$,
repeated application gives
\[
 X^{a_0}z_{a_0,b}=Pz_{0,b},\qquad
 P=\prod_{t=1}^{a_0}C_t\ne0.
\]
Indeed, $a_0<r$, so every factor is nonzero by
Lemma~\ref{lem:firstzero}. Empty products are understood to be $1$.
There are no terms with $a>a_0$, and therefore
\[
 X^{a_0}w_0=P\sum_{b=0}^{b_0}\lambda_{a_0,b}z_{0,b}\in N.
\]
Next, $\beta z_{0,b}=D_{0,b}z_{0,b-1}$. Consequently
$\beta^{b_0}$ kills every term with $b<b_0$, using $D_{0,0}=0$,
whereas
\[
 \beta^{b_0}z_{0,b_0}=Qz_{0,0},\qquad
 Q=\prod_{t=1}^{b_0}D_{0,t}\ne0.
\]
The nonvanishing follows from $b_0<s$ and
Lemma~\ref{lem:firstzero}. Thus
\[
 \beta^{b_0}X^{a_0}w_0
   =\lambda_{a_0,b_0}PQ\,z_{0,0}\in N.
\]
All three scalar factors are nonzero, so $z_{0,0}=v\in N$.
The vector $v$ generates the induced module: explicitly,
$z_{a,b}=\alpha^aY^bv$ for all $0\le a,b<p$. Hence $N=Z$.
We have proved that there is no $H$-submodule strictly between
$W$ and $Z$. Since $W$ is proper, it is a maximal proper submodule.
Consequently, $L(i,j)$ is simple.
\end{proof}
The basis $\{v_{a,b}|0\<a<r,0\<b<s\}$ is called standard basis of $L(i,j)$.
\begin{proposition}\label{prop:explicit}
On the rectangular basis of $L(i,j)$, the actions of $K,\gamma,\alpha,Y,X,\beta$ are \eqref{eq:zweights}--\eqref{eq:Dab}, with $z_{a,b}$ replaced by $v_{a,b}$ and every vector outside the rectangle set equal to zero. The original generators $E,F$ act by
\begin{align}
 Ev_{a,b}&=C_av_{a-1,b}
 -\frac{q^{2j+i-4a}}dD_{a,b}v_{a,b-1},\label{eq:explicitE}\\
 Fv_{a,b}&=v_{a,b+1}+\frac{q^{2-i+2a+2b}}d v_{a+1,b}.
 \label{eq:explicitF}
\end{align}
In particular, $\dim L(i,j)=r(i,j)s(i,j)$.
\end{proposition}
\begin{proof}
The actions descend from the induced $H$-module, so they satisfy every defining relation, including at the boundary of the rectangle. Formula \eqref{eq:recover} gives the actions of $E,F$; here $\gamma^2K$ acts on $v_{a,b}$ by $q^{2j+i-4a}$.
\end{proof}

\begin{lemma}\label{lem:unique}
The induced module $Z(i,j)$ has a unique maximal proper submodule, namely $W(i,j)$, and hence a unique simple quotient.
\end{lemma}
\begin{proof}
Write $Z=Z(i,j)$ and $W=W(i,j)$. Since $\alpha Y=Y\alpha$
and $\alpha^p=Y^p=0$, the assignments $u\mapsto\alpha$ and
$w\mapsto Y$ make $Z$ a module over
\[
 R=\kk[u,w]/(u^p,w^p).
\]
The map $\Phi:R\to Z$, $f\mapsto f(\alpha,Y)v$, is $R$-linear
and sends the basis $\{u^aw^b:0\le a,b<p\}$ of $R$ to the basis
$\{z_{a,b}:0\le a,b<p\}$ of $Z$. Thus $\Phi$ is an isomorphism
from the regular $R$-module onto $Z$.

Let $\mathfrak m=(u,w)\subset R$. Every monomial of total degree
at least $2p-1$ has either its $u$-exponent or its $w$-exponent
at least $p$, so it is zero in $R$. Hence
$\mathfrak m^{2p-1}=0$. If $f=c+n$ with $c\in\kk^{\times}$ and
$n\in\mathfrak m$, the finite geometric-series identity gives
\[
 f^{-1}=c^{-1}\sum_{t=0}^{2p-2}(-c^{-1}n)^t,
\]
because $n^{2p-1}=0$. In particular, every polynomial with
nonzero constant coefficient is invertible in $R$.

Now let $U$ be a proper $H$-submodule of $Z$. It is stable under
$\alpha$ and $Y$, and hence under every polynomial in these
operators. If $U$ contained a vector $f(\alpha,Y)v$ whose
coefficient of $z_{0,0}=v$ were nonzero, then $f$ would have
nonzero constant coefficient. Applying the polynomial operator
$f^{-1}(\alpha,Y)$ would give
\[
 v=f^{-1}(\alpha,Y)f(\alpha,Y)v\in U.
\]
But $z_{a,b}=\alpha^aY^bv$ for every basis vector of $Z$, so
$v\in U$ would imply $U=Z$, a contradiction. Therefore every
proper $H$-submodule is contained in the fixed vector-space
hyperplane
\[
 Z_+=\Span\{z_{a,b}:(a,b)\ne(0,0)\}.
\]
We use $Z_+$ only as a vector subspace; its stability under $H$
is not required.

Let $N$ be the sum of all proper $H$-submodules of $Z$. A sum of
submodules is a submodule, and each of its elements is a finite
sum of vectors from the summands. Since every summand lies in
$Z_+$, so does $N$. Thus $v\notin N$, and $N$ is proper.
Moreover, $N$ contains every proper submodule by construction.
It follows that $N$ is maximal. If $M$ is any maximal
proper submodule, then $M\subseteq N\subsetneq Z$; maximality
forces $M=N$. Thus $N$ is the unique maximal proper submodule.

By Proposition~\ref{prop:simple}, $W$ is a maximal proper
submodule of $Z$. Hence $W=N$. Finally, the kernel of any
surjective homomorphism from $Z$ onto a simple module is maximal
and therefore equals $W$. The first isomorphism theorem then
identifies every simple quotient with $Z/W=L(i,j)$, proving
uniqueness up to isomorphism.
\end{proof}

\begin{theorem}\label{thm:classification}
With $r(i,j)$ and $s(i,j)$ defined by \eqref{eq:r}--\eqref{eq:s},
there is a simple $H$-module $L(i,j)$ with dimension $r(i,j)s(i,j)$. Every simple $H$-module is isomorphic to one of these modules, and
\[
 L(i,j)\cong L(i',j')\quad\Longleftrightarrow\quad
 (i,j)=(i',j')\text{ in }\Z_{2p}^2.
\]
Thus $H$ has exactly $4p^2$ isomorphism classes of simple modules. Explicit generator actions are given in Proposition~\ref{prop:explicit}.
\end{theorem}
\begin{proof}
Existence, simplicity, the generator actions, and the dimension formula follow from Propositions~\ref{prop:simple}--\ref{prop:explicit}.

Let $M$ be any simple $H$-module. By Lemma~\ref{lem:highest}, it contains a nonzero vector $v$ of some highest weight $(i,j)$. Simplicity gives $Hv=M$. The universal property of induction yields a surjective $H$-homomorphism
\[
 Z(i,j)\longrightarrow M,\qquad h\otimes1\longmapsto hv.
\]
Its kernel is maximal. Lemma~\ref{lem:unique} identifies the kernel with $W(i,j)$, proving $M\cong L(i,j)$.

For uniqueness of parameters, we first show that
\begin{equation}
 L(i,j)^X\cap L(i,j)^\beta=\kk v_{0,0}.
 \label{eq:highestline}
\end{equation}
For $1\le a<r$, the coefficient $C_a$ is nonzero, and distinct basis vectors with $a>0$ have distinct images under $X$. Therefore
\[
 L(i,j)^X=\Span\{v_{0,b}:0\le b<s\}.
\]
Within this space, the coefficients $D_{0,b}$ are nonzero for $1\le b<s$, so intersecting with the kernel of $\beta$ gives precisely \eqref{eq:highestline}.

 Let $T$: $L(i,j)\to L(i',j')$ be an $H$-module isomorphism, then $T$ preserves the line in \eqref{eq:highestline} and $T$ commutes with $K,\gamma$. It follows that $q^i=q^{i'}$ and $q^j=q^{j'}$. Since $q$ has order $2p$, the parameters are equal in $\Z_{2p}^2$. The converse is immediate. This proves the classification and the count $4p^2$.
\end{proof}

\section{Indecomposable projective modules}\label{sec:projectives}
We retain the notation of Section~\ref{se:4} and put
\[
 \Lambda=\mathbb Z_{2p}^{\,2},\qquad L_\lambda=L(i,j),\qquad N=\dim H=4p^6,
\]
\[
 d_\lambda=r(i,j)s(i,j)\quad(\lambda=(i,j)),\qquad J=J(H).
\]
All modules in this section are finite dimensional. We construct the
projective covers by a finite linear system and an explicit polynomial.
 Fix an ordering $m_1,\ldots,m_N$ of the PBW basis
\eqref{eq:pbw}. For each $\mu$, order the rectangular basis of $L_\mu$
with $v_{0,0}^\mu$ first. Write
\[
 \kappa_\mu:L_\mu\longrightarrow\kk^{d_\mu}
\]
for the coordinate isomorphism, and let $\mathbf e_0=(1,0,\ldots,0)^{\mathsf T}$.
Thus $\kappa_\mu(v_{0,0}^\mu)=\mathbf e_0$. Define the matrix representation by
\[
 \rho_\mu(h)\kappa_\mu(v)=\kappa_\mu(h\cdot v).
\]
In particular, a matrix $\rho_\mu(h)$ acts on coordinate columns, whereas
$h$ acts on vectors of the abstract module. Let
$E_{00}^\mu=\operatorname{diag}(1,0,\ldots,0)$.

For fixed $\lambda$, solve for scalars $c_{\lambda,t}$ in
\begin{equation}\label{eq:projective-system}
 \sum_{t=1}^{N}c_{\lambda,t}\rho_\mu(m_t)
 =\begin{cases}E_{00}^\lambda,&\mu=\lambda,\\0,&\mu\ne\lambda,
 \end{cases}\qquad\mu\in\Lambda.
\end{equation}
All coefficient matrices are known from Proposition~\ref{prop:explicit}.
Choose any solution, for instance by row reduction with free variables
zero, and define
\begin{align}
 b_\lambda&=\sum_{t=1}^{N}c_{\lambda,t}m_t,\label{eq:projective-lift}\\
 \Phi_N(t)&=\sum_{k=N}^{2N-1}\binom{2N-1}{k}
                         t^k(1-t)^{2N-1-k},\label{eq:idempotent-polynomial}\\
 e_\lambda&=\Phi_N(b_\lambda),\qquad P_\lambda=He_\lambda.
 \label{eq:projective-definition}
\end{align}
We establish the construction and its properties through the following lemmas, and state the classification theorem at the end of this section.

\begin{lemma}\label{lem:semisimple-action-detailed}
The joint action map
\[
 \Psi:H\longrightarrow\prod_{\mu\in\Lambda}M_{d_\mu}(\kk),\qquad
 h\longmapsto(\rho_\mu(h))_\mu
\]
is surjective with kernel $J$. It induces an algebra isomorphism
\[
 \Theta:H/J\xrightarrow{\ \sim\ }\prod_{\mu\in\Lambda}M_{d_\mu}(\kk).
\]
Consequently \eqref{eq:projective-system} has a solution.
\end{lemma}
\begin{proof}
First $J$ annihilates every simple module. Indeed, for a nonzero vector
$v$ of a simple module $L$, the map $H\to L$, $h\mapsto h\cdot v$,
is surjective, so its kernel is a maximal left ideal and contains $J$.
Thus
\[
 v\ne0\Longrightarrow J\subseteq\ker(h\mapsto h\cdot v)
 \Longrightarrow Jv=0,\qquad
 JL=0\Longrightarrow J\subseteq\bigcap_\mu\ker\rho_\mu=\ker\Psi.
\]
Conversely, if $h\in\ker\Psi$, Theorem~\ref{thm:classification} implies
that $h$ annihilates every simple module. For any maximal left ideal $M$,
$H/M$ is simple, so $h\cdot(1+M)=h+M=0$. Consequently
\[
 h\in\ker\Psi
 \Longrightarrow h\cdot(1+M)=0\quad(\forall M\text{ maximal left ideal})
 \Longrightarrow h\in\bigcap_M M=J.
\]
Together these inclusions give
\[
 \ker\Psi=\bigcap_\mu\operatorname{Ann}_H(L_\mu)=J,
 \qquad \Theta(h+J)=(\rho_\mu(h))_\mu.
\]

For surjectivity, the Wedderburn--Artin theorem gives
$H/J\cong\prod_{t=1}^{m}M_{n_t}(D_t)$ with finite-dimensional division
algebras $D_t$ over $\kk$. One can check that each $D_t$ equals $\kk$.
A matrix algebra $M_n(\kk)$ has a unique simple-module isomorphism
class, its natural $n$-dimensional column module. In a product of matrix
algebras, the central identities of the factors give a decomposition
of any module; on a simple module exactly one of them acts as the
identity and the others as zero. Thus the factors of $H/J$ correspond
bijectively to the simple modules $L_\mu$, and their sizes are $d_\mu$.
On each such module the action of its factor identifies that factor
with the full endomorphism algebra. Changing to the fixed rectangular
basis merely conjugates the matrices. In particular,
\[
 H/J\cong\prod_{t=1}^{m}M_{n_t}(\kk),\qquad
 \{(n_t,\text{natural simple module})\}_t
 =\{(d_\mu,L_\mu)\}_\mu
\]
up to reindexing and isomorphism. Therefore
\[
 \Theta:H/J\xrightarrow{\sim}\prod_\mu M_{d_\mu}(\kk),
 \qquad\operatorname{Im}\Psi=\prod_\mu M_{d_\mu}(\kk).
\]
Finally, fix $\lambda\in\Lambda$ and consider the tuple
\[
 T_\lambda=(T_{\lambda,\mu})_{\mu\in\Lambda},\qquad
 T_{\lambda,\mu}=
 \begin{cases}
 E_{00}^\lambda,&\mu=\lambda,\\
 0,&\mu\ne\lambda.
 \end{cases}
\]
Since $\Psi$ is surjective, there exists $b_\lambda\in H$ such that
$\Psi(b_\lambda)=T_\lambda$. Equivalently,
\[
 \rho_\mu(b_\lambda)=
 \begin{cases}
 E_{00}^\lambda,&\mu=\lambda,\\
 0,&\mu\ne\lambda,
 \end{cases}
 \qquad\mu\in\Lambda.
\]
Expand $b_\lambda$ in the fixed PBW basis:
\[
 b_\lambda=\sum_{t=1}^{N}c_{\lambda,t}m_t,
 \qquad c_{\lambda,t}\in\kk.
\]
For every $\mu\in\Lambda$, linearity of $\rho_\mu$ gives
\[
 \sum_{t=1}^{N}c_{\lambda,t}\rho_\mu(m_t)
 =\rho_\mu\!\left(\sum_{t=1}^{N}c_{\lambda,t}m_t\right)
 =\rho_\mu(b_\lambda)
 =\begin{cases}
 E_{00}^\lambda,&\mu=\lambda,\\
 0,&\mu\ne\lambda.
 \end{cases}
\]
Thus equation (23) has a solution.
\end{proof}

\begin{lemma}\label{lem:idempotent-lift-detailed}
The element $e_\lambda$ is idempotent and primitive, and
\begin{equation}\label{eq:projective-idempotent-action}
 \rho_\mu(e_\lambda)=
 \begin{cases}E_{00}^\lambda,&\mu=\lambda,\\0,&\mu\ne\lambda.
 \end{cases}
\end{equation}
In particular $e_\lambda\cdot v_{0,0}^\lambda=v_{0,0}^\lambda$.
\end{lemma}
\begin{proof}
By the proof of lemma \ref{lem:semisimple-action-detailed}, one can check $b_\lambda^2-b_\lambda\in J$.
The descending sequence $H\supset J\supset J^2\supset\cdots$ terminates
in zero by degree $N$: if $J^a=J^{a+1}$, Nakayama's lemma applied to
$J^a$ implies $J^a=0$, while each strict inclusion lowers dimension.
Thus $J^N=0$, and $(b_\lambda^2-b_\lambda)^N=0$.
Every summand of $\Phi_N(t)$ is divisible by $t^N$. Also
\[
 1-\Phi_N(t)=\sum_{k=0}^{N-1}\binom{2N-1}{k}
                 t^k(1-t)^{2N-1-k}
\]
is divisible by $(1-t)^N$. Therefore
$\Phi_N(t)^2-\Phi_N(t)$ is divisible by both $t^N$ and $(1-t)^N$,
and hence by their product, since these polynomials are relatively prime.
More explicitly, for some $Q(t)\in\kk[t]$,
\begin{align*}
 \Phi_N(t)^2-\Phi_N(t)&=t^N(1-t)^NQ(t),\\
 e_\lambda^2-e_\lambda
 &=b_\lambda^N(1-b_\lambda)^NQ(b_\lambda)\\
 &=(-1)^N(b_\lambda^2-b_\lambda)^NQ(b_\lambda)=0.
\end{align*}
Since $\Phi_N(0)=0$ and $\Phi_N(1)=1$, applying each matrix
representation gives \eqref{eq:projective-idempotent-action}.
In coordinates this says
\[
 \kappa_\lambda(e_\lambda\cdot v_{0,0}^\lambda)
 =E_{00}^\lambda\mathbf e_0=\mathbf e_0
 =\kappa_\lambda(v_{0,0}^\lambda),
\]
Hence $e_\lambda\cdot v_{0,0}^\lambda=v_{0,0}^\lambda$.

Write $\overline e_\lambda=e_\lambda+J$. By lemma \ref{lem:semisimple-action-detailed} and (\ref{eq:projective-idempotent-action}), $\Theta(\overline e_\lambda)=E_{00}^\lambda$, so
$\overline e_\lambda\ne0$. To check primitivity, suppose
\[
 \overline e_\lambda=u+v,\qquad u^2=u,\quad v^2=v,\quad uv=vu=0.
\]
Then $\overline e_\lambda u=u$ and $\overline e_\lambda v=v$.
Thus all components of $u,v$ outside the $\lambda$-factor vanish.
Writing $u_\lambda,v_\lambda$ for their matrices in that factor gives
\[
 E_{00}^\lambda u_\lambda=u_\lambda,\qquad
 E_{00}^\lambda v_\lambda=v_\lambda,
 \qquad
 \operatorname{Im}u_\lambda,\operatorname{Im}v_\lambda
 \subseteq\operatorname{Im}E_{00}^\lambda,
\] which imply $ u_\lambda=0$ or $ v_\lambda=0$, and hence
 $u=0$ or $v=0$. Thus $\overline e_\lambda$ is primitive.

Finally, if $e_\lambda=f+g$ with $f,g$ orthogonal idempotents in $H$,
then $\overline e_\lambda=\overline f+\overline g$ is an orthogonal
idempotent decomposition in $H/J$. By primitivity, one summand,
say $\overline f$, is zero. Consequently
\[
 \overline f=0\Longrightarrow f\in J
 \Longrightarrow f=f^N\in J^N=0.
\]
The argument for $\overline g=0$ is identical. Thus $e_\lambda$ is primitive.
\end{proof}

\begin{lemma}\label{lem:projective-ideal-detailed}
Let $P_\lambda=He_\lambda$. Then  $P_\lambda$ is an indecomposable projective module, and
\[
 JP_\lambda=Je_\lambda=He_\lambda\cap J.
\]
\end{lemma}
\begin{proof}
Since $e_\lambda$ is idempotent and primitive, $P_\lambda$ is an indecomposable projective module.

Since $J$ is a right ideal and $1\in H$, $JH=J$.
Hence
\begin{align*}
 JP_\lambda=J(He_\lambda)&=(JH)e_\lambda=Je_\lambda.\\
\end{align*}
For the intersection, the two inclusions are
\begin{align*}
 je_\lambda\in Je_\lambda&\Longrightarrow je_\lambda\in He_\lambda\ \text{and}\ je_\lambda\in J
          \Longrightarrow Je_\lambda\subseteq He_\lambda\cap J,\\
 z\in He_\lambda\cap J&\Longrightarrow z=he_\lambda,
 \quad ze_\lambda=h{e_\lambda}^2=he_\lambda=z\\
 &\Longrightarrow z=ze_\lambda\in Je_\lambda.
\end{align*}
Therefore $JP_\lambda=Je_\lambda=He_\lambda\cap J$.
\end{proof}

\begin{lemma}\label{lem:top-projective-detailed}
There is an $H$-module isomorphism
\[
 P_\lambda/JP_\lambda\cong L_\lambda.
\]
The map
\[
 \pi_\lambda:P_\lambda\longrightarrow L_\lambda,\qquad
 he_\lambda\longmapsto h\cdot v_{0,0}^\lambda
\]
is a projective cover with kernel $JP_\lambda$.
\end{lemma}
\begin{proof}
The restriction of the quotient map $H\to H/J$ gives the surjection
\[
 He_\lambda\longrightarrow(H/J)\overline e_\lambda,\qquad
 he_\lambda\longmapsto(h+J)\overline e_\lambda.
\]
Its kernel is $He_\lambda\cap J=JP_\lambda$, by  lemma \ref{lem:projective-ideal-detailed}.
Hence
\[
 P_\lambda/JP_\lambda
 \cong(H/J)\overline e_\lambda
 \cong M_{d_\lambda}(\kk)E_{00}^\lambda
 \cong L_\lambda.
\]
For clarity, we describe the other two maps explicitly.
Under $\Theta$, the element $(h+J)\overline e_\lambda$ has only one
nonzero component, $\rho_\lambda(h)E_{00}^\lambda$ in position $\lambda$.
Taking this component is the second isomorphism. Give the matrix
left ideal the $H$-action $h\cdot C=\rho_\lambda(h)C$.
For the third isomorphism use
\[
 g:M_{d_\lambda}(\kk)E_{00}^\lambda\longrightarrow L_\lambda,
 \qquad C\longmapsto\kappa_\lambda^{-1}(C\mathbf e_0).
\]
Writing the ordered module basis as
$(v_0^\lambda,v_1^\lambda,\ldots,v_{d_\lambda-1}^\lambda)$,
where $v_0^\lambda=v_{0,0}^\lambda$, we have
\[
 C=(c_{ab})\in M_{d_\lambda}(\kk)E_{00}^\lambda
 \Longleftrightarrow c_{ab}=0\quad(b\ne0),
 \qquad
 g(C)=\sum_{a=0}^{d_\lambda-1}c_{a0}v_a^\lambda.
\]
Its inverse is explicitly
\[
 g^{-1}\Bigl(\sum_ac_av_a^\lambda\Bigr)
   =\sum_ac_aE_{a0}^\lambda.
\]
Thus $g$ is bijective. To check $H$-linearity, observe first that
\[
 C=AE_{00}^\lambda
 \Longrightarrow\rho_\lambda(h)C
   =(\rho_\lambda(h)A)E_{00}^\lambda
   \in M_{d_\lambda}(\kk)E_{00}^\lambda.
\]
Then
\[
 \kappa_\lambda(g(h\cdot C))
 =\rho_\lambda(h)C\mathbf e_0
 =\kappa_\lambda(h\cdot g(C)).
\] Hence $g$ is an $H$-module isomorphism.

The composite sends $he_\lambda+JP_\lambda$ to
\[
 \kappa_\lambda^{-1}(\rho_\lambda(h)E_{00}^\lambda\mathbf e_0)
 =\kappa_\lambda^{-1}(\rho_\lambda(h)\mathbf e_0)
 =h\cdot v_{0,0}^\lambda.
\]
Consequently $\pi_\lambda$ is well defined, surjective, and has kernel
$JP_\lambda$. One can also check well-definedness directly: if
$(h-h')e_\lambda=0$, then
$(h-h')\cdot v_{0,0}^\lambda=(h-h')e_\lambda\cdot v_{0,0}^\lambda=0$.

Finally, if $U+JP_\lambda=P_\lambda$, put $M=P_\lambda/U$.
Then
\[
 M={(U+JP_\lambda)}/U=JM=J^2M=\cdots=J^NM=0
 \Longrightarrow U=P_\lambda.
\]
This is exactly the assertion that $JP_\lambda$ is superfluous:
adding it to a proper submodule never gives the whole module.
Since $P_\lambda$ is projective, the surjection with this kernel
is a projective cover.
\end{proof}

\begin{theorem}\label{thm:projective-covers}
The $P_\lambda$ give all pairwise nonisomorphic indecomposable projective
modules. Moreover,
\begin{equation}\label{eq:regular-projective-decomposition}
 {}_HH\cong\bigoplus_{\lambda\in\Lambda}P_\lambda^{\,d_\lambda}.
\end{equation}
\end{theorem}
\begin{proof}
Since $P_{\lambda}$ are projective cover of $L_{\lambda}$, $P_{\lambda}$ are pairwise nonisomorphic indecomposable projective modules.  Hence we write
$H\cong\bigoplus_\lambda P_\lambda^{n_\lambda}$.
Taking quotients by $J$ gives
\[
 H/J\cong\bigoplus_\lambda(P_\lambda/JP_\lambda)^{n_\lambda}
      \cong\bigoplus_\lambda L_\lambda^{n_\lambda}.
\]
On the other hand, the left regular module of $M_d(\kk)$ is the
direct sum of its $d$ column left ideals $M_d(\kk)E_{aa}$, each
isomorphic to the natural $d$-dimensional simple module.
Lemma~\ref{lem:semisimple-action-detailed} consequently gives
$H/J\cong\bigoplus_\lambda L_\lambda^{d_\lambda}$.
Comparing simple multiplicities proves $n_\lambda=d_\lambda$.
\end{proof}

We now determine which of the simple modules are projective and compute
the dimensions of their projective covers. Write $P_{i,j}=P_\lambda$
for $\lambda=(i,j)$. The following calculation applies to both parities
of $p$.

Put $B^-=\langle\alpha,Y,K,\gamma\rangle$. By
Lemma~\ref{lem:pbw}, this subalgebra has basis
\[
 \{\alpha^aY^bK^k\gamma^\ell:
       0\le a,b<p,\ k,\ell\in\Z_{2p}\},
\]
and $H$ is free as a left $B^-$-module with basis
$\{X^a\beta^b:0\le a,b<p\}$. For $(u,v)\in\Lambda$, let
$\kk_{u,v}^{\,\mathrm r}$ be the one-dimensional right $B^-$-module
on which $\alpha,Y$ act as zero and $K,\gamma$ act by $q^u,q^v$.
Define the right induced module and its left dual by
\[
 Z^{\mathrm r}(u,v)=\kk_{u,v}^{\,\mathrm r}\otimes_{B^-}H,
 \qquad Q(u,v)=\operatorname{Hom}_{\kk}(Z^{\mathrm r}(u,v),\kk).
\]
The left action on this dual is $(hf)(z)=f(zh)$. No antipode is
used in this passage from right modules to left modules.

\begin{lemma}\label{lem:projective-dimension-filtration}
The right regular module $H_H$ has a filtration in which every
$Z^{\mathrm r}(u,v)$ occurs $p^2$ times. Moreover, for every
$(i,j)\in\Lambda$,
\[
 \dim_{\kk}P_{i,j}
 =p^2\sum_{(u,v)\in\Lambda}[Q(u,v):L(i,j)],
\]
where brackets denote composition multiplicity.
\end{lemma}
\begin{proof}
Set $
 F^dB^-=\Span_{\kk}\{\alpha^aY^bK^k\gamma^\ell:
 0\le a,b<p,\ k,\ell\in\Z_{2p},\ a+b\ge d\}.
$
The defining relations give
\[
 (F^dB^-)\alpha+(F^dB^-)Y\subseteq F^{d+1}B^-,
 \qquad
 (F^dB^-)K+(F^dB^-)\gamma\subseteq F^dB^-.
\]
Thus these spaces are right ideals, with $F^0B^-=B^-$ and
$F^{2p-1}B^-=0$. For $a+b=d$, put
\[
 w_{a,b;u,v}=\sum_{k,\ell=0}^{2p-1}q^{-uk-v\ell}
                 \overline{\alpha^aY^bK^k\gamma^\ell}
 \quad ((u,v)\in\Lambda),
\]
where the bar denotes the class in $F^dB^-/F^{d+1}B^-$. Shifting
$k$ or $\ell$ in the sum gives
\[
 w_{a,b;u,v}K=q^u w_{a,b;u,v},\qquad
 w_{a,b;u,v}\gamma=q^v w_{a,b;u,v},\qquad
 w_{a,b;u,v}\alpha=w_{a,b;u,v}Y=0.
\]
For fixed $(a,b)$, these $(2p)^2$ nonzero vectors have distinct
joint eigenvalues and form a basis of the corresponding PBW span.
Hence, as right $B^-$-modules,
\[
 F^dB^-/F^{d+1}B^-
 \cong \bigoplus_{\substack{0\le a,b<p\\a+b=d}}
       \ \bigoplus_{(u,v)\in\Lambda}\kk_{u,v}^{\,\mathrm r}.
\]
Refinement gives a composition series of $B^-_{B^-}$ with $p^2$
copies of each $\kk_{u,v}^{\,\mathrm r}$, one for each pair $(a,b)$.
Since $H$ is free as a left $B^-$-module, $-\otimes_{B^-}H$
is exact. Under $B^-\otimes_{B^-}H\cong H_H$, the induced filtration
has subquotients
\[
 \kk_{u,v}^{\,\mathrm r}\otimes_{B^-}H=Z^{\mathrm r}(u,v),
\]
each occurring $p^2$ times.

Fix $(i,j)$ and put $e=e_{i,j}$. For any right $H$-module map
$f:N\to N'$, the identities $xe=x$ for $x\in Ne$ and
$f(ne)=f(n)e$ yield
\[
 \ker(f|_{Ne})=(\ker f)e,\qquad f(Ne)=f(N)e.
\]
Thus $N\mapsto Ne$ is exact as a functor to vector spaces.
By \eqref{eq:projective-idempotent-action}, $e_{i,j}$ has dimension
one on the right simple module $L(i,j)^*$ and dimension zero on the
other right simples. Here $L(i,j)^*$ is the ordinary linear dual,
with $(fh)(z)=f(hz)$. Hence
\[
 \dim_{\kk}Ne=[N:L(i,j)^*],
\]
where brackets denote composition multiplicity. Applying this
to the preceding filtration, and then taking linear duals, yields
\begin{equation}\label{eq:projective-dimension-multiplicity}
 \dim_{\kk}P_{i,j}=\dim_{\kk}He
 =p^2\sum_{(u,v)\in\Lambda}[Q(u,v):L(i,j)].
\end{equation}

\end{proof}

\begin{lemma}\label{lem:dual-induced-rectangle-filtration}
Put $r_0=r(u,v)$ and $s_0=s(u,v)$. Partition $\{0,\ldots,p-1\}$
into $[0,r_0-1]$ and, when $r_0<p$, $[r_0,p-1]$ for the first
index, and analogously at $s_0$ for the second index.
Then $Q(u,v)$ has a filtration with one simple subquotient for
each resulting rectangle. The rectangle with lower corner
$(a_0,b_0)$ gives the simple module
\[
 L(u-2a_0-2b_0,\ v-a_0+b_0).
\]
In particular, $Q(u,v)$ has at most four composition factors,
counted with multiplicity.
\end{lemma}
\begin{proof}
The right module
$Z^{\mathrm r}(u,v)$ has basis $z^{\mathrm r}_{a,b}
=1\otimes X^a\beta^b$, $0\le a,b<p$. Let $\xi_{a,b}$ be its
dual basis in $Q(u,v)$. Define
\[
 C_a(u,v)=q^{-2}(a)_{q^{-2}}
             (q^{2v+u-2(a-1)}-1),\qquad
 D_{0,b}(u,v)=(b)_{q^2}(q^{-u}-q^{-2v-2(b-1)}).
\]
These are the coefficients of \eqref{eq:Ca}--\eqref{eq:Dab}
with parameters $(u,v)$. Reordering products in the right induced
module gives
\[
\begin{aligned}
 z^{\mathrm r}_{a,b}K&=q^{u-2a-2b}z^{\mathrm r}_{a,b},
 &z^{\mathrm r}_{a,b}\gamma&=q^{v-a+b}z^{\mathrm r}_{a,b},\\
 z^{\mathrm r}_{a,b}X&=q^{-2b}z^{\mathrm r}_{a+1,b},
 &z^{\mathrm r}_{a,b}\beta&=z^{\mathrm r}_{a,b+1},\\
 z^{\mathrm r}_{a,b}\alpha&=q^{2b}C_a(u,v)z^{\mathrm r}_{a-1,b},
 &z^{\mathrm r}_{a,b}Y&=q^{2a}D_{0,b}(u,v)z^{\mathrm r}_{a,b-1}.
\end{aligned}
\]
For the $\alpha$ formula, commute $\alpha$ past $\beta^b$
and then use the expansion of $X^a\alpha$ from
$X\alpha=q^{-2}\alpha X+q^{-2}(\gamma^2K-1)$; the term with
$\alpha$ on the far left vanishes in the induced module.
For the $Y$ formula, iterate
$\beta Y=Y\beta+K^{-1}-\gamma^{-2}$.
Moving its correction term to the left of $X^a$ contributes
$q^{2a}$, since
\[
 X^aK^{-1}=q^{2a}K^{-1}X^a,\qquad
 X^a\gamma^{-2}=q^{2a}\gamma^{-2}X^a.
\]
The term with $Y$ on the far left vanishes. Dualizing gives
\begin{equation}\label{eq:dual-right-induced-actions}
\begin{aligned}
 K\xi_{a,b}&=q^{u-2a-2b}\xi_{a,b},
 &\gamma\xi_{a,b}&=q^{v-a+b}\xi_{a,b},\\
 X\xi_{a,b}&=q^{-2b}\xi_{a-1,b},
 &\beta\xi_{a,b}&=\xi_{a,b-1},\\
 \alpha\xi_{a,b}&=q^{2b}C_{a+1}(u,v)\xi_{a+1,b},
 &Y\xi_{a,b}&=q^{2a}D_{0,b+1}(u,v)\xi_{a,b+1}.
\end{aligned}
\end{equation}
Every vector with an index outside $\{0,\ldots,p-1\}$ is zero.

Write $r_0=r(u,v)$ and $s_0=s(u,v)$. The internal zeros of the
raising coefficients in \eqref{eq:dual-right-induced-actions}
occur exactly at the cuts $a=r_0$ and $b=s_0$, whenever these
numbers are less than $p$. Split the first index interval into
$[0,r_0-1]$ and, if $r_0<p$, $[r_0,p-1]$; split the second in
the same way at $s_0$. These intervals give at most four
rectangles.

Let \(\mathcal R\) be the set of these rectangles. For \(R\in\mathcal R\)
write \(R=[a_0,a_1]\times[b_0,b_1]\) and set
\[
 V_R=\operatorname{span}_{\kk}\{\xi_{a,b}:(a,b)\in R\}.
\]
Order the rectangles by their lower corners, using the coordinatewise
partial order. By \eqref{eq:dual-right-induced-actions}, we have
\[
 gV_R\subseteq\sum_{S\le R}V_S
 \qquad(g\in\{K,\gamma,X,\beta,\alpha,Y\}).
\]
Indeed, \(K,\gamma,\alpha,Y\) preserve \(V_R\), since the raising
coefficients vanish at the cuts, and \(X,\beta\) can only move a basis
vector to the same rectangle or to a smaller one. Choose a listing
\(R_1,\ldots,R_m\) such that \(R_i<R_j\) implies \(i<j\). Then
\[
 0=F_0\subset F_1\subset\cdots\subset F_m=Q(u,v),
 \qquad F_j=\bigoplus_{i=1}^{j}V_{R_i},
\]
is a filtration by \(H\)-submodules, with the direct sums taken as vector spaces.

Fix \(j\), write \(R_j=[a_0,a_1]\times[b_0,b_1]\), and let
\(\bar\xi_{a,b}\) be the image of \(\xi_{a,b}\) in \(F_j/F_{j-1}\)
for \((a,b)\in R_j\). These images form a basis. The actions are
given by \eqref{eq:dual-right-induced-actions}, with any term leaving
\(R_j\) interpreted as zero in the quotient.

Let \(U\) be a nonzero submodule of \(F_j/F_{j-1}\) and take
\(0\ne w=\sum_{(a,b)\in R_j}c_{a,b}\bar\xi_{a,b}\in U\).
Put
\[
 a_*=\max\{a:c_{a,b}\ne0\text{ for some }b\},
 \qquad b_*=\max\{b:c_{a_*,b}\ne0\}.
\]
Terms with \(a<a_*\) vanish under \(X^{a_*-a_0}\), so
\begin{align*}
 X^{a_*-a_0}w
 &=\sum_{b=b_0}^{b_1}c_{a_*,b}q^{-2b(a_*-a_0)}\bar\xi_{a_0,b},\\
 \beta^{b_*-b_0}X^{a_*-a_0}w
 &=c_{a_*,b_*}q^{-2b_*(a_*-a_0)}\bar\xi_{a_0,b_0}.
\end{align*}
The last coefficient is nonzero, hence \(\bar\xi_{a_0,b_0}\in U\).
For every \((a,b)\in R_j\), the raising coefficients along the
corresponding steps are nonzero, and
\[
 Y^{b-b_0}\alpha^{a-a_0}\bar\xi_{a_0,b_0}
 =
 \left(\prod_{t=a_0+1}^{a}q^{2b_0}C_t(u,v)\right)
 \left(\prod_{t=b_0+1}^{b}q^{2a}D_{0,t}(u,v)\right)
 \bar\xi_{a,b},
\]
where an empty product is \(1\). Thus \(U=F_j/F_{j-1}\).

Finally, with \(\eta=\bar\xi_{a_0,b_0}\), the action formulas give
\[
 X\eta=\beta\eta=0,\qquad
 K\eta=q^{u-2a_0-2b_0}\eta,\qquad
 \gamma\eta=q^{v-a_0+b_0}\eta.
\]
Theorem~\ref{thm:classification} therefore gives
\[
 F_j/F_{j-1}\cong L(u-2a_0-2b_0,\ v-a_0+b_0).
\]
Since \(m=|\mathcal R|\le4\), this filtration has at most four
simple subquotients, counted with multiplicity.

\end{proof}

\begin{lemma}\label{lem:projective-dimension-formula}
For every $(i,j)\in\Lambda$, put $r=r(i,j)$ and $s=s(i,j)$.
Then
\[
 \dim_{\kk}P_{i,j}=p^2\varepsilon(r)\varepsilon(s),\qquad
 \varepsilon(t)=
 \begin{cases}
 1,&t=p,\\
 2,&1\le t<p.
 \end{cases}
\]
\end{lemma}
\begin{proof}
By Lemmas~\ref{lem:projective-dimension-filtration}
and~\ref{lem:dual-induced-rectangle-filtration}, it remains to
count all rectangles giving the fixed simple module $L(i,j)$.
Fix now a target $L(i,j)$, and put $r=r(i,j)$, $s=s(i,j)$.
If $i$ is odd, the source parameter $u$ must be odd as well;
then $r_0=s_0=p$ and there is only one rectangle. It contributes
precisely for $(u,v)=(i,j)$.

If $i$ is even, a rectangle of width $r$ has first-coordinate
origin
\[
 a_0\in A_r:=
 \begin{cases}
 \{0\},&r=p,\\
 \{0,p-r\},&r<p,
 \end{cases}
\]
according as it is the lower or upper interval. Similarly its
second-coordinate origin belongs to $A_s$. For each pair
$(a_0,b_0)\in A_r\times A_s$, the highest weight of the rectangle
forces, uniquely in $\Z_{2p}^2$,
\[
 u=i+2a_0+2b_0,\qquad v=j+a_0-b_0.
\]
These sources really have the required cuts. Indeed,
\[
\begin{aligned}
 r(u,v)-1&\equiv r-1+2a_0\pmod p,\\
 s(u,v)-1&\equiv s-1+2b_0\pmod p.
\end{aligned}
\]
For $a_0=0$ this gives $r_0=r$; for $a_0=p-r>0$ it gives
$r_0=a_0$, so the upper interval has length $p-r_0=r$.
The argument for the second index is identical. Counting
rectangles with multiplicity, even if two of them occur in the
same source module, therefore proves
\begin{equation}\label{eq:projective-dimension-explicit}
 \dim_{\kk}P_{i,j}=p^2\varepsilon(r)\varepsilon(s),
 \qquad
 \varepsilon(t)=
 \begin{cases}
 1,&t=p,\\
 2,&1\le t<p.
 \end{cases}
\end{equation}
For odd $i$, this formula also holds, since $r=s=p$.

\end{proof}

\begin{theorem}\label{thm:simple-projective-classification}
The simple module $L(i,j)$ is projective if and only if
$r(i,j)=s(i,j)=p$. Equivalently, the simple projective modules are:
\begin{enumerate}[label=\textup{(\alph*)}]
 \item all $L(i,j)$ with $i$ odd;
 \item all $L(i,j)$ with $i$ even and
 \[
 j+i/2\equiv-1\pmod p,\qquad i/2-j\equiv-1\pmod p.
 \]
\end{enumerate}
Every such module has dimension $p^2$ and is isomorphic to its
projective cover. There are exactly $2p^2+2\gcd(2,p)$
isomorphism classes of simple projective modules.
\end{theorem}
\begin{proof}
If $r=s=p$, \eqref{eq:projective-dimension-explicit} gives
$\dim_{\kk}P_{i,j}=p^2=\dim_{\kk}L(i,j)$, so the projective-cover
surjection is an isomorphism. Conversely, if either $r<p$ or
$s<p$, then
\[
 \dim_{\kk}L(i,j)=rs<p^2\le\dim_{\kk}P_{i,j}.
\]
A projective simple module is isomorphic to its projective cover:
the cover splits by projectivity, and its indecomposable domain
forces the kernel to vanish. Thus $L(i,j)$ is not projective
in this case.

For odd $i$ both edge lengths are $p$. For even $i$, the two
congruences in the statement are precisely the conditions
$r=s=p$ in \eqref{eq:r}--\eqref{eq:s}. There are $2p^2$ choices
with $i$ odd. For $i$ even, write $t=i/2\in\Z_p$. The two
congruences are equivalent to
\[
 2(t+1)\equiv0\pmod p,\qquad j\equiv-1-t\pmod p.
\]
There are $\gcd(2,p)$ choices of $t$, and each admits two lifts
of $j$ to $\Z_{2p}$. The count follows from the uniqueness of
parameters in Theorem~\ref{thm:classification}.
\end{proof}

\begin{theorem}\label{thm:nonsimple-projective-dimensions}
Every nonsimple finite-dimensional indecomposable projective
$H$-module is isomorphic to a unique $P_{i,j}$ with $i$ even and
$(r(i,j),s(i,j))\ne(p,p)$. Put $g=\gcd(2,p)$. Their dimensions
and numbers of isomorphism classes are
\[
\begin{array}{c|c|c}
 \text{condition on }r,s&\dim_{\kk}P_{i,j}&
       \text{number of classes}\\\hline
 r=p,\ s<p\ \text{or }r<p,\ s=p&2p^2&4p-4g\\
 r<p,\ s<p&4p^2&2p^2-4p+2g
\end{array}
\]
where $r=r(i,j)$ and $s=s(i,j)$.
\end{theorem}
\begin{proof}
By Theorem~\ref{thm:projective-covers}, the modules $P_{i,j}$
exhaust the indecomposable projectives, without repetition.
Theorem~\ref{thm:simple-projective-classification} identifies the
simple ones. The dimensions in the remaining cases follow from
\eqref{eq:projective-dimension-explicit}.

Among the $2p^2$ parameters with $i$ even, exactly $2p$ satisfy
$r=p$: each $t=i/2\in\Z_p$ fixes $j\equiv-1-t\pmod p$, with
two lifts in $\Z_{2p}$. Similarly, exactly $2p$ satisfy $s=p$.
Their intersection has $2g$ elements by the preceding theorem.
Thus exactly
\[
 (2p-2g)+(2p-2g)=4p-4g
\]
parameters have exactly one edge length equal to $p$.
The number with both edge lengths less than $p$ is
\[
 2p^2-(2p+2p-2g)=2p^2-4p+2g,
\]
as asserted.
\end{proof}

\section{Strict wildness}\label{sec:wildness}
Let $R=\kk\langle x,y\rangle$. We write $R$-mod and $H$-mod for the
categories of finite-dimensional left modules. An $R$-mod
is a triple $M=(V;A,B)$, where $A,B\in\operatorname{End}_{\kk}(V)$
are the actions of $x,y$. Since $R$ is free on $x$, $y$, any pair of endomorphisms defined such a module.  In a basis of $V$ these are arbitrary square matrices.
We construct a fully faithful exact functor
\[F:R\text{-mod}\longrightarrow H\text{-mod}\]
and represent it by a bimodule free of finite rank on the right.
The main theorem is stated after the construction.

By Section 3, one know that
the presentation of $H$ are given by  Lemmas \ref{2.2}, \ref{lem:modified}
 and \[\quad K^{2p}=1,\quad K\gamma=\gamma K,\quad
 K\alpha=q^{-2}\alpha K,\quad K\beta=q^2\beta K.
\]

\begin{lemma}\label{lem:free-pair-hom}
For $M=(V;A,B)$ and $M'=(V';A',B')$,
\[
 \operatorname{Hom}_R(M,M')
 =\{P\in\operatorname{Hom}_{\kk}(V,V'):PA=A'P,\ PB=B'P\}.
\]
\end{lemma}
\begin{proof} It is easy to see that $\operatorname{Hom}_R(M,M')$ is determined by matrix $P \in \operatorname{Hom}_{\kk}(V,V')$.
Let $f_{P}\in \operatorname{Hom}_R(M,M')$ then $f_P(rz)=rf_P(z)$. Taking $r=x,y$ gives
\[P(Az)=A'Pz,\qquad P(Bz)=B'Pz.\]
Conversely, these identities imply,
by induction on the length of a word $w$ in $x,y$,
\[
 f_P(w\cdot z)=w\cdot f_P(z).
\]
If $\dim V=n$ and $\dim V'=m$, then $P$ is an $m\times n$ matrix.
It is an isomorphism precisely when $m=n$, $P$ is invertible, and
$A'=PAP^{-1}$, $B'=PBP^{-1}$.
\end{proof}

\begin{lemma}\label{lem:wild-simple-endpoints}
Suppose $p\ge3$. The simple modules
\[
 U_1=L(-2,-1),\qquad U_2=L(-2,p-1),\qquad U_3=L(-2,1)
\]
are pairwise nonisomorphic and have rectangular sizes
$(p-1,1),(p-1,1),(1,p-1)$, respectively. With bases
$u_a,v_a$ and $w_b$, indexed by $0\le a,b\le p-2$, they satisfy
\begin{gather*}
 Xu_0=Xv_0=0,\qquad \alpha u_{p-2}=\alpha v_{p-2}=0,\\
 \beta w_0=0,\qquad Yw_{p-2}=0,\\
 YU_1=\beta U_1=YU_2=\beta U_2=0,\qquad \alpha U_3=XU_3=0.
\end{gather*}
Their endpoint weights are
\[
\begin{array}{c|cc}
 &K&\gamma\\\hline
 u_0&q^{-2}&q^{-1}\\
 u_{p-2}&q^2&-q\\
 v_0&q^{-2}&-q^{-1}\\
 v_{p-2}&q^2&q\\
 w_0&q^{-2}&q\\
 w_{p-2}&q^2&-q^{-1}
\end{array}
\]
\end{lemma}
\begin{proof}
All parameters are taken modulo $2p$. For $i=-2$ we may use $i/2=-1$
in the residue formulas. Thus the pairs $(r,s)$ are
\[
 (1+\res{-2},1+\res{0}),\quad
 (1+\res{p-2},1+\res{-p}),\quad
 (1+\res{0},1+\res{-2}),
\]
which give the asserted sizes. The three second parameters are distinct
modulo $2p$ for $p\ge3$, so the simples are nonisomorphic by
Theorem~\ref{thm:classification}.
On $U_1,U_2$ the second index is always zero; hence
$Yv_{a,0}=v_{a,1}=0$ and $\beta v_{a,0}=0$, using $D_{a,0}=0$.
At their first endpoint $Xv_{0,0}=0$ because $C_0=0$, and at the
other endpoint $\alpha v_{p-2,0}=v_{p-1,0}=0$ by the quotient boundary.
On $U_3$ the first index is zero; consequently
$\alpha v_{0,b}=v_{1,b}=0$ and $Xv_{0,b}=0$.
Its endpoints satisfy $\beta v_{0,0}=0$ and
$Yv_{0,p-2}=v_{0,p-1}=0$.
Finally substitute these indices into
$Kv_{a,b}=q^{i-2a-2b}v_{a,b}$ and
$\gamma v_{a,b}=q^{j-a+b}v_{a,b}$, using $q^p=-1$,
to obtain the displayed weight table.
\end{proof}

For $p\ge3$ and $M=(V;A,B)$ define the vector space
\begin{equation}\label{eq:wild-space-detailed}
 F(M)=V_+\oplus V_-\oplus\bigoplus_{i=1}^{3}(U_i\otimes_{\kk}V).
\end{equation}
Here $V_+,V_-$ are copies of $V$, with vectors $t_+(z),t_-(z)$.
On each $U_i\otimes V$ define $h(u\otimes z)=(hu)\otimes z$.
This is a multiplicity-space construction: $V$ has no specified
$H$-action, and no coproduct is used here.
On the two other summands set
\[
 Kt_\pm(z)=t_\pm(z),\qquad\gamma t_\pm(z)=\pm t_\pm(z),
\]
and prescribe
\begin{equation}\label{eq:wild-action-table}
\begin{array}{c|cc}
 &t_+(z)&t_-(z)\\\hline
 \alpha&u_0\otimes z&v_0\otimes Az\\
 X&v_{p-2}\otimes z&u_{p-2}\otimes z\\
 Y&w_0\otimes z&0\\
 \beta&0&w_{p-2}\otimes Bz
\end{array}
\end{equation}
The direct sum in \eqref{eq:wild-space-detailed} is initially only a
vector-space decomposition. In particular $V_+,V_-$ need not be submodules.

\begin{proposition}\label{prop:wild-module-valid}
For arbitrary $A,B$, these formulas define an $H$-module $F(M)$.
\end{proposition}
\begin{proof}
Put $N_V=\bigoplus_{i=1}^3(U_i\otimes V)$. The relations already hold on $N_V$,
so it suffices to check $t_\pm(z)$ using
above relations.
For a vector of joint weight $(k,g)$, the required target weights
under the four operators are
\[
\begin{array}{c|cc}
 &K&\gamma\\\hline
 \alpha&q^{-2}k&q^{-1}g\\
 X&q^2k&qg\\
 Y&q^{-2}k&qg\\
 \beta&q^2k&q^{-1}g
\end{array}
\]
Every nonzero entry in \eqref{eq:wild-action-table} has precisely
this weight by Lemma~\ref{lem:wild-simple-endpoints}. For instance,
\begin{align*}
 K\alpha t_+(z)&=q^{-2}u_0\otimes z=q^{-2}\alpha Kt_+(z),\\
 \gamma\alpha t_-(z)&=-q^{-1}v_0\otimes Az
                       =q^{-1}\alpha\gamma t_-(z).
\end{align*}
This verifies all Cartan weight relations; their finite-order and
commutativity relations are immediate on the two summands.
Moreover
\[
 (\gamma^2K-1)t_\pm(z)=0,
 \qquad (K^{-1}-\gamma^{-2})t_\pm(z)=0.
\]
For the first corrected relation we compute
\begin{align*}
 X\alpha t_+(z)&=Xu_0\otimes z=0,&
 \alpha Xt_+(z)&=\alpha v_{p-2}\otimes z=0,\\
 X\alpha t_-(z)&=Xv_0\otimes Az=0,&
 \alpha Xt_-(z)&=\alpha u_{p-2}\otimes z=0.
\end{align*}
Hence $X\alpha-q^{-2}\alpha X=q^{-2}(\gamma^2K-1)$ holds there.
For the second, the four expressions are
\begin{align*}
 \beta Yt_+(z)&=\beta w_0\otimes z=0,&Y\beta t_+(z)&=0,\\
 \beta Yt_-(z)&=0,&Y\beta t_-(z)&=Yw_{p-2}\otimes Bz=0.
\end{align*}
Thus $\beta Y-Y\beta=K^{-1}-\gamma^{-2}$ holds.
Every one of $Y\alpha,\alpha Y,XY,YX,X\beta,\beta X,
\alpha\beta,\beta\alpha$ annihilates both top summands, by the
zero actions on $U_i$. For example,
$XYt_+(z)=Xw_0\otimes z=0$ and
$YXt_+(z)=Yv_{p-2}\otimes z=0$; the other products are checked
from the same table. This verifies all remaining mixed relations.
For $g\in\{\alpha,X,Y,\beta\}$ we have
\[
 g(V_+\oplus V_-)\subseteq N_V,\qquad
 g^{p-1}N_V=0\quad\Longrightarrow\quad g^pF(M)=0.
\]
All defining relations have now been verified.
\end{proof}

\begin{lemma}\label{lem:wild-hom-blocks}
For $p\ge3$, every $H$-homomorphism $T:F(M)\to F(M')$ has the form
\[
 T=P_+\oplus P_-\oplus\bigoplus_{i=1}^{3}(1_{U_i}\otimes C_i)
\]
for linear maps $P_+,P_-,C_i:V\to V'$.
\end{lemma}
\begin{proof}
Write
\[
 N_V=\bigoplus_{i=1}^3(U_i\otimes_{\kk}V),\qquad
 N_{V'}=\bigoplus_{i=1}^3(U_i\otimes_{\kk}V').
\]
Thus $F(M)=V_+\oplus V_-\oplus N_V$, and similarly for $F(M')$.
The spaces $V_\pm$ are copies of $V$, identified with $V$ by
$t_\pm:V\to V_\pm$; write $t'_\pm:V'\to V'_\pm$ for the corresponding
identifications in the target. On each $U_i\otimes_{\kk}V$ the action
is $h(u\otimes z)=(hu)\otimes z$. In particular, $N_V$ is an
$H$-submodule, although $V_+$ and $V_-$ need not be $H$-submodules.

First we separate these summands using the actions of $K$ and
$\gamma$. Denote their action operators on $F(M)$ by $K_M,\gamma_M$,
and those on $F(M')$ by $K_{M'},\gamma_{M'}$. Since $T$ is $H$-linear,
\[
 TK_M=K_{M'}T,\qquad T\gamma_M=\gamma_{M'}T.
\]
For $c\in\kk$, let
\[
 F(M)_c=\{\xi\in F(M):K_M\xi=c\xi\}
\]
be the $K$-eigenspace with eigenvalue $c$, and define $F(M')_c$
in the same way. If $\xi\in F(M)_c$, then
\[
 K_{M'}T(\xi)=T(K_M\xi)=cT(\xi).
\]
Consequently $T(F(M)_c)\subseteq F(M')_c$.

The displayed bases of the $U_i$ consist of $K$-eigenvectors, with
eigenvalues in
\[
 S=\{q^{-2a}:1\le a\le p-1\}.
\]
Since $q^2$ has order $p$, we have $1\notin S$. On the other hand,
$K$ acts as the identity on $V_+\oplus V_-$. Hence
\[
 F(M)_1=V_+\oplus V_-,
 \qquad N_V=\bigoplus_{c\in S}F(M)_c,
\]
with the analogous equalities for $M'$. It follows that
\[
 T(V_+\oplus V_-)\subseteq V'_+\oplus V'_-,
 \qquad T(N_V)\subseteq N_{V'}.
\]
Within the first of these spaces, $\gamma$ distinguishes the two
summands:
\[
 \gamma_Mt_+(z)=t_+(z),\qquad
 \gamma_Mt_-(z)=-t_-(z).
\]
Its eigenvalues $1$ and $-1$ are distinct over $\kk$. Commutation
with $\gamma$ therefore gives $T(V_\pm)\subseteq V'_\pm$.
There are thus unique $\kk$-linear maps $P_\pm:V\to V'$ such that
\[
 Tt_\pm(z)=t'_\pm(P_\pm z)\qquad(z\in V).
\]

It remains to determine the restriction $T|_{N_V}$. Choose bases
$e_1,\ldots,e_n$ of $V$ and $f_1,\ldots,f_m$ of $V'$, where
$n=\dim_{\kk}V$ and $m=\dim_{\kk}V'$. For each $i$ and $a$, there
are unique linear maps $T^{j,i}_{ba}:U_i\to U_j$ satisfying
\[
 T(u\otimes e_a)=
 \sum_{j=1}^3\sum_{b=1}^m T^{j,i}_{ba}(u)\otimes f_b
 \qquad(u\in U_i).
\]
Here $i,a$ specify the source copy of $U_i$, and $j,b$ specify the
target copy of $U_j$. Equivalently, $T^{j,i}_{ba}$ is obtained by
including $U_i$ as $U_i\otimes\kk e_a$, applying $T|_{N_V}$, and
projecting onto $U_j\otimes\kk f_b$, identified with $U_j$.
These inclusions and projections are $H$-linear, because the action
on every tensor summand is confined to its first factor.

Explicitly, for $h\in H$ and $u\in U_i$, $H$-linearity yields
\begin{align*}
 \sum_{j,b}T^{j,i}_{ba}(hu)\otimes f_b
 &=T((hu)\otimes e_a)\\
 &=hT(u\otimes e_a)
 =\sum_{j,b}hT^{j,i}_{ba}(u)\otimes f_b.
\end{align*}
Uniqueness of the expansion in the target direct sum implies
\[
 T^{j,i}_{ba}(hu)=hT^{j,i}_{ba}(u),
 \qquad
 T^{j,i}_{ba}\in\operatorname{Hom}_H(U_i,U_j).
\]
By Lemma~\ref{lem:wild-simple-endpoints}, the $U_i$ are pairwise
nonisomorphic simple modules. By Schur's lemma, we have
\[
 \operatorname{Hom}_H(U_i,U_j)=0\quad(i\ne j),
 \qquad
 \operatorname{End}_H(U_i)=\kk\,1_{U_i},
\]
where $1_{U_i}$ denotes the identity map of $U_i$.
Thus $T^{j,i}_{ba}=0$ for $j\ne i$, and there are scalars
$c^{(i)}_{ba}\in\kk$ such that
\[
 T^{i,i}_{ba}=c^{(i)}_{ba}1_{U_i}.
\]
Define $C_i:V\to V'$ by
\[
 C_i(e_a)=\sum_{b=1}^m c^{(i)}_{ba}f_b.
\]
The matrix $(c^{(i)}_{ba})$ is therefore the $m\times n$ matrix of
$C_i$ in the chosen bases. Substitution gives
\[
 T(u\otimes e_a)
 =\sum_b c^{(i)}_{ba}u\otimes f_b
 =u\otimes C_i(e_a).
\]
By linearity, for every $u\in U_i$ and $z\in V$,
\[
 T(u\otimes z)=u\otimes C_i(z),
 \qquad
 T|_{U_i\otimes V}=1_{U_i}\otimes C_i.
\]
Combining these restrictions proves
\[
 T=P_+\oplus P_-\oplus
       \bigoplus_{i=1}^3(1_{U_i}\otimes C_i).
\]
This is a decomposition of a linear map with respect to the stated
vector-space direct sums. Further equalities among $P_+,P_-$ and
the $C_i$ require the connecting actions of $\alpha,X,Y,\beta$;
they are established in the next proposition.
\end{proof}

\begin{proposition}\label{prop:wild-full-faithful}
For $p\ge3$, every $H$-homomorphism $T:F(M)\to F(M')$ is uniquely
of the form
\[
 F(P)=P\oplus P\oplus\bigoplus_{i=1}^{3}(1_{U_i}\otimes P),
 \qquad PA=A'P,\quad PB=B'P.
\]
These maps define a fully faithful exact functor $F:R\text{-mod}\to H\text{-mod}$.
\end{proposition}
\begin{proof}
Start with the five maps in Lemma~\ref{lem:wild-hom-blocks}.
Commuting with $\alpha$ on $t_+(z)$ gives
\[
 u_0\otimes C_1z=T\alpha t_+(z)
               =\alpha Tt_+(z)=u_0\otimes P_+z,
 \qquad C_1=P_+.
\]
Here $u_0\ne0$, so $u_0\otimes-$ is injective (apply a linear
functional taking $u_0$ to $1$). The same computation with $X,Y$
on $t_+(z)$ yields
\[
 v_{p-2}\otimes C_2z=v_{p-2}\otimes P_+z,\qquad
 w_0\otimes C_3z=w_0\otimes P_+z,
 \qquad C_2=C_3=P_+.
\]
Commuting with $X$ on $t_-(z)$ similarly gives
$u_{p-2}\otimes C_1z=u_{p-2}\otimes P_-z$, hence
\[
 P_+=P_-=C_1=C_2=C_3=:P.
\]
The two parameter-dependent entries now give
\begin{align*}
 v_0\otimes PAz&=T\alpha t_-(z)=\alpha Tt_-(z)=v_0\otimes A'Pz,\\
 w_{p-2}\otimes PBz&=T\beta t_-(z)=\beta Tt_-(z)=w_{p-2}\otimes B'Pz.
\end{align*}
Thus $PA=A'P$ and $PB=B'P$.
Conversely, these two identities make the displayed diagonal map
commute with every entry of \eqref{eq:wild-action-table}; it also
commutes with $K,\gamma$ and with the actions on the $U_i$ factors.
It is therefore $H$-linear, proving the natural bijection
\begin{equation}\label{eq:wild-hom}
 \operatorname{Hom}_H(F(M),F(M'))
 \cong\{P:PA=A'P,\ PB=B'P\}
 =\operatorname{Hom}_R(M,M').
\end{equation}
Moreover $F(1_M)=1_{F(M)}$ and $F(QP)=F(Q)F(P)$ by the diagonal
formula. Hence $F$ is a fully faithful functor.
For a short exact sequence in $R$-mod, applying $F$ gives the direct
sum of two copies of the underlying sequence and its tensor products
with $U_1,U_2,U_3$ over the field $\kk$. These are exact sequences of
vector spaces, and all maps are $H$-linear. Thus the sequence is
exact in $H$-mod, proving exactness of $F$.
\end{proof}

\begin{proposition}\label{prop:wild-p-two}
For $p=2$ there is a fully faithful exact functor $F:R\text{-mod}\to H\text{-mod}$
with $\dim F(M)=3\dim M$.
\end{proposition}
\begin{proof}
When $p=2$, we have $q^2=-1$ and $q^{-1}=-q$. The preceding
construction requires a modification, since $U_2=U_3$ in this case.
Let $M=(V;A,B)$ be a finite-dimensional left $R$-module, so that
$A,B\in\operatorname{End}_{\kk}(V)$ give the actions of $x,y$,
respectively. Set
\[
 F(M)=V_0\oplus V_-\oplus V_+,
\]
where the three summands are copies of $V$. Denote their fixed
linear identifications with $V$ by
\[
 t:V\longrightarrow V_0,\qquad
 s_-:V\longrightarrow V_-,\qquad
 s_+:V\longrightarrow V_+.
\]
In particular, $t(z),s_-(z),s_+(z)$ denote vectors in distinct
summands, not scalar multiples of a common vector.

Define $K$ and $\gamma$ by their scalar actions
\[
\begin{array}{c|ccc}
 &V_0&V_-&V_+\\\hline
 K&1&-1&-1\\
 \gamma&1&q^{-1}&q
\end{array}
\]
and define the other four generators by
\[
\begin{aligned}
 \alpha t(z)&=s_-(z),& \beta t(z)&=s_-(Az),\\
 Xt(z)&=s_+(z),& Yt(z)&=s_+(Bz).
\end{aligned}
\]
Each of $\alpha,\beta,X,Y$ is defined to be zero on
$W_V:=V_-\oplus V_+$. These prescriptions extend linearly to
$F(M)$, whose dimension is $3\dim_{\kk}V$.

We first verify that the prescribed operators satisfy the defining
relations of $H$. The operators $K,\gamma$ are diagonal with respect
to the displayed decomposition, so
\[
 K^4=\gamma^4=1,\qquad K\gamma=\gamma K.
\]
Here and below, $1$ denotes the identity operator on $F(M)$.
The three joint weights, meaning the pairs of eigenvalues for
$K,\gamma$, are
\[
 (1,1)\text{ on }V_0,\qquad
 (-1,q^{-1})\text{ on }V_-,\qquad
 (-1,q)\text{ on }V_+.
\]
The weight relations for $\alpha,\beta$ require a vector of weight
$(k,g)$ to be sent to weight $(-k,q^{-1}g)$, whereas those for
$X,Y$ require weight $(-k,qg)$. Thus the prescribed maps from
$V_0$ into $V_-$ or $V_+$ have exactly the required weights.
For example,
\[
 K\alpha t(z)=-s_-(z)=q^{-2}\alpha Kt(z),
 \qquad
 \gamma Yt(z)=q\,s_+(Bz)=qY\gamma t(z).
\]
The same check applies to $\beta$ and $X$, and all weight relations
hold trivially on $W_V$, where the four operators vanish.
This verifies the relations involving one of $K,\gamma$ and one
of $\alpha,\beta,X,Y$.

For any $a,b\in\{\alpha,\beta,X,Y\}$ we have
\[
 bF(M)\subseteq W_V,\qquad aW_V=0,
 \qquad\text{hence }ab=0.
\]
In particular,
\[
 \alpha^2=\beta^2=X^2=Y^2=0,
\]
and the relations
\[
 Y\alpha=\alpha Y,\quad XY=YX,\quad
 X\beta=q^2\beta X,\quad \alpha\beta=q^{-2}\beta\alpha
\]
hold because both sides vanish. To verify the two relations with
correction terms, the scalar actions above give
\[
 \gamma^2K=1,\qquad K^{-1}=\gamma^{-2}.
\]
Indeed, the eigenvalues of $\gamma^2K$ on the three summands are
$1,q^{-2}(-1),q^2(-1)$, all equal to $1$; those of both
$K^{-1}$ and $\gamma^{-2}$ are $1,-1,-1$.
Consequently both sides of
\[
 X\alpha-q^{-2}\alpha X=q^{-2}(\gamma^2K-1),
 \qquad
 \beta Y-Y\beta=K^{-1}-\gamma^{-2}
\]
vanish. All defining relations are satisfied, and $F(M)$ is
therefore an $H$-module. This verification places no restriction
on $A$ and $B$.

Next let $M'=(V';A',B')$ and let $T:F(M)\to F(M')$ be
$H$-linear. Use $t',s'_-,s'_+$ for the corresponding identifications
in $F(M')$. Since $T$ commutes with the actions of $K$ and $\gamma$,
it preserves joint weights: if $K\xi=k\xi$ and $\gamma\xi=g\xi$,
then
\[
 K T(\xi)=T(K\xi)=kT(\xi),\qquad
 \gamma T(\xi)=T(\gamma\xi)=gT(\xi).
\]
The three joint weights are distinct, since $1\ne-1$ and
$q\ne q^{-1}$. Hence there are unique linear maps
$P_0,P_-,P_+:V\to V'$ with
\[
 Tt(z)=t'(P_0z),\qquad
 Ts_-(z)=s'_-(P_-z),\qquad
 Ts_+(z)=s'_+(P_+z).
\]
Commuting with $\alpha$ and $X$ on $t(z)$ gives
\begin{align*}
 s'_-(P_-z)&=T\alpha t(z)=\alpha Tt(z)=s'_-(P_0z),\\
 s'_+(P_+z)&=TXt(z)=XTt(z)=s'_+(P_0z).
\end{align*}
Since $s'_-$ and $s'_+$ are injective, it follows that
\[
 P_-=P_0=P_+=:P.
\]
The remaining two generators then give
\begin{align*}
 s'_-(PAz)&=T\beta t(z)=\beta Tt(z)=s'_-(A'Pz),\\
 s'_+(PBz)&=TYt(z)=YTt(z)=s'_+(B'Pz).
\end{align*}
Thus $PA=A'P$ and $PB=B'P$, which are precisely the conditions for
$P:M\to M'$ to be $R$-linear.

Conversely, if $P:V\to V'$ satisfies these two identities, define
$F(P)$ by
\[
 F(P)t(z)=t'(Pz),\qquad
 F(P)s_\pm(z)=s'_\pm(Pz).
\]
This map preserves the three weight spaces and hence commutes with
$K,\gamma$. The preceding computations, read in reverse, show that
it commutes with $\alpha,X,\beta,Y$ on $V_0$. It also commutes
with them on $W_V$, since they annihilate $W_V$ and $W_{V'}$.
Therefore $F(P)$ is $H$-linear. We obtain a bijection
\[
 \operatorname{Hom}_R(M,M')
 \longrightarrow\operatorname{Hom}_H(F(M),F(M')),
 \qquad P\longmapsto F(P).
\]
The definitions give
\[
 F(1_M)=1_{F(M)},\qquad F(QP)=F(Q)F(P)
\]
for composable $R$-linear maps $P,Q$. Thus $F$ is a fully faithful
functor.

Finally, consider a short exact sequence of finite-dimensional
$R$-modules
\[
 0\longrightarrow M_1\xrightarrow{P}M_2
   \xrightarrow{Q}M_3\longrightarrow0.
\]
Under the fixed identifications of the three copies, $F(P)$ and
$F(Q)$ are $P\oplus P\oplus P$ and $Q\oplus Q\oplus Q$.
The resulting sequence of vector spaces is the direct sum of three
copies of the underlying short exact sequence, and is therefore
exact. All its maps are $H$-linear, as proved above, so
\[
 0\longrightarrow F(M_1)\xrightarrow{F(P)}F(M_2)
   \xrightarrow{F(Q)}F(M_3)\longrightarrow0
\]
is exact in $H$-mod. This proves exactness and completes the proof.
\end{proof}

\begin{lemma}\label{lem:wild-bimodule-and-consequences}
In both cases $F\cong\mathcal M\otimes_R-$ for an $H$--$R$-bimodule
$\mathcal M$ free of finite rank on the right. Moreover $F$ reflects
isomorphism and preserves indecomposability.
\end{lemma}
\begin{proof}
We first construct the representing bimodule. For $p\ge3$, take
\[
 \mathcal M=R_+\oplus R_-\oplus
       \bigoplus_{i=1}^3(U_i\otimes_{\kk}R),
\]
where $R_+,R_-$ are copies of the right regular module $R_R$,
with elements denoted by $t_+(r),t_-(r)$. Give this space the
componentwise right action
\[
 t_\pm(r)s=t_\pm(rs),\qquad
 (u\otimes r)s=u\otimes rs\qquad(r,s\in R).
\]
Define the left actions of the generators of $H$ by the same
formulas as in the construction of $F(M)$, replacing $V$ by $R$
and $A,B$ by the operators
\[
 L_x(r)=xr,\qquad L_y(r)=yr.
\]
In particular,
\[
 \alpha t_-(r)=v_0\otimes xr,\qquad
 \beta t_-(r)=w_{p-2}\otimes yr.
\]
On $U_i\otimes_{\kk}R$ the left action is
$h(u\otimes r)=(hu)\otimes r$.

For $p=2$, take instead
\[
 \mathcal M=R_0\oplus R_-\oplus R_+,
\]
with the componentwise right regular action, and use the formulas
of Proposition~\ref{prop:wild-p-two} with $A=L_x$ and $B=L_y$.
Thus, for example,
\[
 \alpha t(r)=s_-(r),\quad \beta t(r)=s_-(xr),\quad
 Xt(r)=s_+(r),\quad Yt(r)=s_+(yr).
\]
In both cases all prescribed operators are right $R$-linear.
For instance, in the first construction,
\[
 \alpha\bigl(t_-(r)s\bigr)
 =v_0\otimes x(rs)
 =v_0\otimes(xr)s
 =\bigl(\alpha t_-(r)\bigr)s.
\]
The same calculation applies to $L_y$; the other action maps use
only scalar multiplication, identity maps, zero maps, or operators
on the $U_i$ factor, and therefore also commute with the right action.
The defining relations of $H$ hold by the calculations in
Propositions~\ref{prop:wild-module-valid} and~\ref{prop:wild-p-two}.
 We have consequently constructed an $H$--$R$-bimodule:
\[
 h(ms)=(hm)s
 \qquad(h\in H,\ m\in\mathcal M,\ s\in R).
\]

Choosing a $\kk$-basis of each $U_i$ identifies
$U_i\otimes_{\kk}R$ with $R^{p-1}$ as right $R$-modules. Hence
\[
 \mathcal M_R\cong
 \begin{cases}
 R^{\,2+3(p-1)}=R^{\,3p-1},&p\ge3,\\
 R^{\,3},&p=2.
 \end{cases}
\]
Here the exponent denotes the number of copies of the right regular
module; in particular, $\mathcal M$ is free of finite rank on the
right, although it is infinite-dimensional over $\kk$.

We next identify its tensor functor with $F$. Let $M=(V;A,B)$.
For $p\ge3$, define
\[
 \theta_M:\mathcal M\otimes_R M\longrightarrow F(M)
\]
on the summands by
\[
 \theta_M(t_\pm(r)\otimes z)=t_\pm(rz),\qquad
 \theta_M((u\otimes r)\otimes z)=u\otimes rz.
\]
Here $rz$ denotes the given left $R$-action on $M$, so that
$xz=Az$ and $yz=Bz$. These formulas are balanced over $R$:
for example, $(rs)z=r(sz)$ ensures that
$(u\otimes rs)\otimes z$ and $(u\otimes r)\otimes sz$
have the same image. They therefore define a linear map on the
tensor product. Its inverse is given by
\[
 t_\pm(z)\longmapsto t_\pm(1)\otimes z,\qquad
 u\otimes z\longmapsto(u\otimes1)\otimes z.
\]
The tensor product has left $H$-action
$h(m\otimes z)=(hm)\otimes z$, well defined by the bimodule
identity above. The map $\theta_M$ intertwines this action with
the prescribed action on $F(M)$. For the parameter-dependent
$\alpha$-action, for example,
\begin{align*}
 \theta_M\bigl(\alpha(t_-(r)\otimes z)\bigr)
 &=\theta_M((v_0\otimes xr)\otimes z)\\
 &=v_0\otimes(xr)z
 =v_0\otimes A(rz)
 =\alpha\,\theta_M(t_-(r)\otimes z).
\end{align*}
The $\beta$-action follows by replacing $x,A$ with $y,B$;
the remaining generators follow directly from their formulas.
For $p=2$, the same construction is given by
\[
 \theta_M(t(r)\otimes z)=t(rz),\qquad
 \theta_M(s_\pm(r)\otimes z)=s_\pm(rz),
\]
with inverse obtained by putting $r=1$; its $H$-linearity is checked
in exactly the same manner.

For an $R$-linear map $P:M\to M'$, the equality $P(rz)=rP(z)$
gives
\[
 \theta_{M'}(1_{\mathcal M}\otimes P)=F(P)\theta_M.
\]
Thus the $\theta_M$ form a natural isomorphism
$F\cong\mathcal M\otimes_R-$. Right freeness also implies that
$\mathcal M\otimes_R-$ is exact, providing another proof of the
exactness of $F$.

It remains to prove the two assertions about modules. Suppose
$S:F(M)\to F(M')$ is an $H$-module isomorphism. By fullness there
are $R$-linear maps $P:M\to M'$ and $Q:M'\to M$ satisfying
\[
 F(P)=S,\qquad F(Q)=S^{-1}.
\]
Functoriality gives
\[
 F(QP)=S^{-1}S=F(1_M),\qquad
 F(PQ)=SS^{-1}=F(1_{M'}).
\]
Faithfulness, namely injectivity on each set of homomorphisms,
therefore implies
\[
 QP=1_M,\qquad PQ=1_{M'}.
\]
Thus $P$ is an isomorphism and $M\cong M'$, proving that $F$
reflects isomorphism.

Finally, suppose $M$ is indecomposable. Full faithfulness induces
an isomorphism of rings
\[
 \operatorname{End}_R(M)\xrightarrow{\ \sim\ }
       \operatorname{End}_H(F(M)),\qquad f\longmapsto F(f).
\]
If $F(M)=N_1\oplus N_2$ with both $N_i$ nonzero, the projection
$e$ onto $N_1$ along $N_2$ is $H$-linear and satisfies
\[
 e^2=e,\qquad e\ne0,\qquad e\ne1.
\]
Let $f$ be its preimage in $\operatorname{End}_R(M)$. Since the
displayed map is a ring isomorphism,
\[
 f^2=f,\qquad f\ne0,\qquad f\ne1.
\]
Both $\operatorname{Im}f$ and $\ker f$ are $R$-submodules. For every
$z\in M$ we have
\[
 z=f(z)+(z-f(z)),\qquad f(z-f(z))=f(z)-f^2(z)=0.
\]
Moreover, if $w=f(u)\in\ker f$, then
$w=f^2(u)=f(w)=0$. Consequently
\[
 M=\operatorname{Im}f\oplus\ker f.
\]
The first summand is nonzero because $f\ne0$. The second is nonzero
because $f\ne1$: choose $z$ with $z-f(z)\ne0$.
This contradicts indecomposability of $M$, so $F(M)$ is
indecomposable.
\end{proof}

\begin{theorem}\label{thm:strict-wildness}
For every integer $p\ge2$, the algebra $H$ is strictly wild.
There is a fully faithful exact representation embedding
$R\text{-mod}\to H\text{-mod}$ represented by a right-free bimodule
of rank $3p-1$ for $p\ge3$ and rank $3$ for $p=2$.
\end{theorem}
\begin{proof}
Propositions~\ref{prop:wild-module-valid}
and~\ref{prop:wild-full-faithful} give the fully faithful exact functor for $p\>3$, and
Proposition~\ref{prop:wild-p-two} gives it for $p=2$. Lemma~\ref{lem:wild-bimodule-and-consequences} identifies the functor with tensoring by a right-free bimodule of the stated rank. It also shows that the functor reflects isomorphism and preserves indecomposability. This gives a strict representation-embedding, so $H$ is strictly  wild and hence wild.
\end{proof}

\centerline{ACKNOWLEDGMENTS}

This work is supported by the National Natural Science Foundation of China  (No.12201545).\\

{\bf Author Contributions:} H. S ,  Y. S, X. L  and H. C wrote the  main manuscript text. All authors reviewed the manuscript.

{\bf Materials availability:} Not applicable.

{\bf Code availability:} Not applicable.

\section*{Declarations}
{\bf Funding:} This work is supported by National Natural Science  Foundation  of China (No.12201545).\\

{\bf Data availability statement:} This manuscript has no associated data.\\

{\bf Ethics approval:} Not applicable.\\

{\bf Conflict of interest:} The authors declare that there is no conflict of interest.

\end{document}